\documentclass[%
 reprint,%
 amssymb, amsmath,%
 aip,cha,%
]{revtex4-1}

\usepackage{color}%
\usepackage{bm}%
\usepackage{graphicx}%
\usepackage{caption}
\usepackage{float}
\usepackage{subfig}
\usepackage{dcolumn}%
\usepackage{bm}%
\usepackage[colorlinks=true,linkcolor=blue]{hyperref}%
\expandafter\ifx\csname package@font\endcsname\relax\else
 \expandafter\expandafter
 \expandafter\usepackage
 \expandafter\expandafter
 \expandafter{\csname package@font\endcsname}%
\fi
\begin{document}

\newcommand\bes{\begin{eqnarray}}
\newcommand\ees{\end{eqnarray}}
\newcommand\bess{\begin{eqnarray*}}
\newcommand\eess{\end{eqnarray*}}
\newcommand{\ve}{\varepsilon}
\newtheorem{definition}{Definition}
\newtheorem{theorem}{Theorem}[section]
\newtheorem{lemma}{Lemma}[section]
\newtheorem{proposition}{Proposition}[section]
\newtheorem{remark}{Remark}[section]
\newtheorem{corollary}{Corollary}[section]
\newtheorem{example}{Example}[section]
\title[Stochastic Predator-Prey Model]
{\bf Stochastic mixed-mode oscillations in a three-species predator-prey model}
\author{Susmita Sadhu}
\author{Christian Kuehn}
\address{Department of Mathematics,
Georgia College \& State University, Milledgeville, GA 31061, USA}
\email{susmita.sadhu@gcsu.edu}
\address{Department of Mathematics,
Technical University of Munich, Germany}
\email{ckuehn@ma.tum.de}


\thispagestyle{empty}


\begin{abstract}
\noindent    
The effect of demographic stochasticity, in the form of Gaussian white noise, 
in a predator-prey model with one fast and two slow variables is studied. 
We derive the stochastic differential equations (SDEs) from a discrete 
model. For suitable parameter values, the deterministic drift part of the model 
admits a folded node singularity and exhibits a singular Hopf bifurcation. We 
focus on the parameter regime near the Hopf bifurcation, where small amplitude 
oscillations exist as stable dynamics in the absence of noise. In this regime, 
the stochastic model admits noise-driven mixed-mode oscillations (MMOs), which capture 
the intermediate dynamics between two cycles of population outbreaks. We perform 
numerical simulations to calculate the distribution of the random number of small 
oscillations between successive spikes for varying noise intensities and distance 
to the Hopf bifurcation. We also study the effect of noise on a suitable 
Poincar\'e map. Finally, we prove that the stochastic model can be transformed 
into a normal form near the folded node, which can be linked to recent results
on the interplay between deterministic and stochastic small amplitude oscillations.
The normal form can also be used to study the parameter influence on the noise
level near folded singularities.
\end{abstract}

\maketitle


{\bf Population dynamics is subject to demographic stochasticity. Dynamical models
are crucial to understand the effect of this noise. We illustrate
in the context of a predator-prey model the entire multi-stage process, starting 
from the derivation of the model using probabilistic fundamentals, to numerical simulations
and concluding with links to abstract mathematical results on normal forms.
We find interesting effects related to oscillation patterns, which can
be constructed from a geometric phase space decomposition. Parts of the dynamics
are slow, and essentially metastable, until sudden fast tipping events occur,
either induced by noise or by deterministic slow drift. The resulting MMOs 
can be viewed as switching in a multi-stable system.}


\section{Introduction}
\label{sec:intro}

Stochasticity or variability plays an important role in understanding the 
dynamics of predator-prey populations. In many cases, it helps to explain 
observable population phenomena such as random fluctuations in population 
densities of small mammals and forest insects~\cite{AC, KBJP, M1}. 
Higher per-capita growth rate of the prey than its predators can be modeled 
to naturally give rise to systems of singularly-perturbed differential 
equations~\cite{BD1, H, RM, S3}. Noise in such models can lead 
to several interesting dynamical effects, which are not anticipated by their 
deterministic counterpart~\cite{S1}. 

There are several ways  to model variability in population ecology~\cite{C1}. 
These include variations between individuals usually referred to as demographic 
stochasticity, or 
variations caused by temporal changes in the environment, usually called 
environmental variations. In this paper, we will assume that variations 
in between individuals, follows a certain 
probability distribution. In this form of modeling, all 
individuals in the same population are treated as identical in all aspects 
of phenotype, and the randomness enters as a within-individual process 
having the same probability distribution for every individual regardless 
of its phenotype. We build a stochastic model, and derive when the population 
sizes are large, that the dynamics is approximated by a system of
SDEs. A similar 
approach to modeling two or more interacting populations with demographic 
stochasticity naturally built-in has been considered in~\cite{A1, A2}.

We consider the stochastic dynamics between three interacting populations, 
two predators competing for their common prey, under two different 
timescales. The deterministic model, i.e., the drift terms
in the SDEs, has been extensively studied in~\cite{S2, S3, SCT1}. Of the many 
interesting dynamics that the deterministic model shows are MMOs, which are 
concatenations of small and large amplitude 
oscillations \cite{BKW, DGKKOW, KuehnBook}. The MMO cycles are ecologically 
significant as they qualitatively represent natural fluctuations in the 
population densities of species as observed in the wild. Taking demographic 
stochasticity into account, we obtain the interactions in the framework of 
a fast-slow system of SDEs. We focus our study regarding noise 
effects in a parameter regime near the singular Hopf bifurcation, where 
noise-induced MMOs are observed. To the best of our knowledge, noise-induced 
MMOs in three interacting populations have not been previously explored. Our
results show that the distance to the singular Hopf bifurcation and noise
level are crucial parameters for formation of different oscillation 
patterns. Large noise levels lead to 
stochastic extinction events, while very low noise levels preclude  
frequent noise-induced large amplitude oscillations (LAOs). Furthermore, 
parts of canard-type orbits near unstable manifolds are kicked away by
small perturbations, thus ``regularizing'' the Poincar\'e map. 
Therefore, we conclude that optimal finite-size noise 
is beneficial in those predator-prey systems, where switching 
between alternating ecological regimes may be desirable.
Hence, our model is a more complex example of the classical paradigms 
of stochastic/coherent resonance 
\cite{BenziSuteraVulpiani,NicolisNicolis,GammaitoniHaenggiJungMarchesoni}, 
which are frequently studied in neuroscience \cite{LindnerSchimansky-Geier1}, and 
also provides a new direction to explore other noise-induced mechanisms  \cite{KuehnNetworks}.

In this paper we make an attempt to understand  the interplay between 
singularities in fast-slow systems  
and noise, which has gained considerable attention 
recently~\cite{BerglundGentz,BerglundGentzKuehn1}. We therefore
transform our system to a suitable normal form locally; this seems to be
the first three-dimensional realistic SODE model with multiplicative noise, where such a transformation
near folded nodes has been carried out. Finally, we compare our model
results to theoretical predictions. In summary, we have provided a paradigmatic
example from applications, where we combine modeling, stochastic fast-slow 
systems, numerical simulations and algebra to gain insight into oscillatory
patterns and multi-stability. 

The paper is organized as follows. In Section~\ref{sec:deterministic}, 
we introduce the deterministic model and review relevant theory of 
fast-slow systems. Existence of folded nodes and canards are also discussed 
in that section. In Section~\ref{sec:main}, we build our main fast-slow stochastic 
model, we
perform numerical simulations and 
discuss the results. In Section~\ref{sec:fold_analysis}, the stochastic model 
is transformed to its normal form and the implications are discussed. The derivation
of the normal form is provided in Appendix~\ref{sec:appnd}.


\section{The Deterministic Model}
\label{sec:deterministic}

Before we can start with the full stochastic model, we have to 
understand the deterministic aspects of the dynamics. This part forms
the skeleton of our analysis as we are interested in 
small stochastic perturbations of deterministic dynamics. 
The deterministic model in its non-dimensional form with the 
derivation details provided in~\cite{SCT1} reads as
\begin{eqnarray}
\label{nondim1}    
\left\{
\begin{array}{lll} \frac{\textnormal{d}x}{\textnormal{d}t}
&= x'&= x\left(1-x-\frac{y}{\beta_1+x}
-\frac{z}{\beta_2+x}\right),\\
   \frac{\textnormal{d}y}{\textnormal{d}t}&=y'&=\zeta_1 
	y\left(\frac{x}{\beta_1+x}-c\right),\\
    \frac{\textnormal{d}z}{\textnormal{d}t}&=z' 
		&= \zeta_2  z\left(\frac{x}{\beta_2+x}-d-hz\right),
       \end{array} 
\right. 
\end{eqnarray}
where $x$ is the normalized population density of the prey and $y$ 
and $z$ are the normalized  population densities of the predators.
We always restrict our study to the first quadrant 
$x\geq 0$, $y\geq 0$, $z\geq 0$ due to the natural non-negativity 
assumption in population dynamics. 
The parameters $\zeta_1$, $\zeta_2$, $\beta_1$, $\beta_2$, $c$, $d$ 
and $h$ represent the following:

\begin{itemize}
 \item $\zeta_1, \ \zeta_2$ denote the ratios of the maximum per 
capita growth rates of the predators $y$ and $z$ respectively 
to the per-capita growth rate of the prey. We will assume 
that  $0<\zeta_1, \zeta_2\ll 1$ and fix $\zeta_1=\zeta_2=:\zeta$.
 \item The parameters $\beta_1$ and $\beta_2$ are dimensionless 
semi-saturation constants measured against the prey's carrying 
capacity. We will assume that the predators  will reach the half 
of their maximum predation rates before 
the prey population reaches its carrying capacity implying 
that $0< \beta_1,\ \beta_2<1$.
 \item The parameters $c$ and $d$ denote the ratios of the growth 
rates of the predators to their death rates. We will make the
natural  assumption
that $0<c,\ d<1$. 
 \item The parameter $h$ is associated with the intraspecific 
competition in the class of $z$. In our context, 
$h>1$ will indicate the presence of  a strong intraspecific 
competition while $0<h<1$ will indicate a weaker competition.
\end{itemize}

\begin{figure}[ht!]%
    \begin{center}
    \subfloat[The critical 
		manifold.]{{\includegraphics[width=5.5cm]{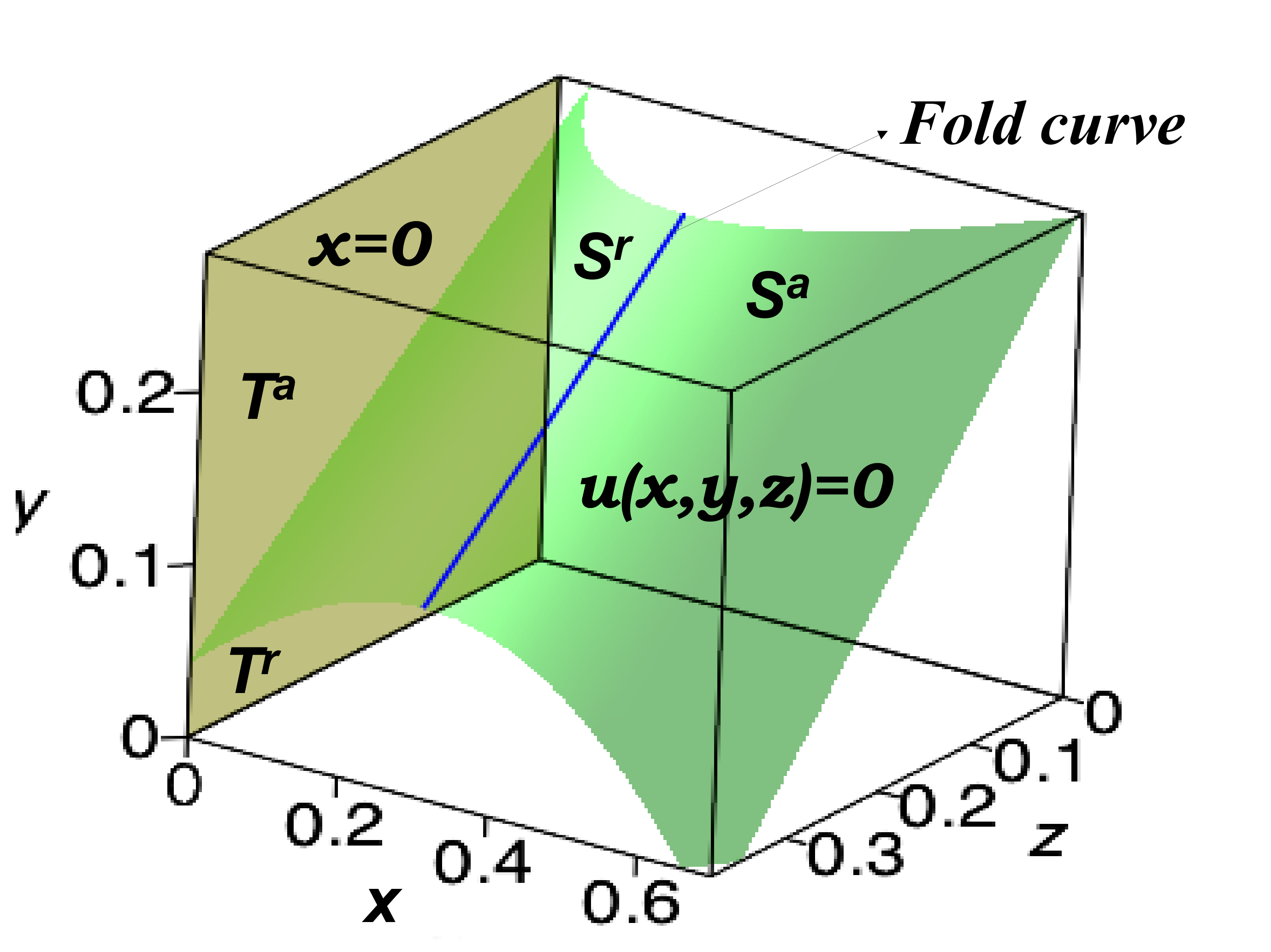} }}%
    \qquad
    \subfloat[Existence of a global 
		return mechanism]{{\includegraphics[width=6.2cm]{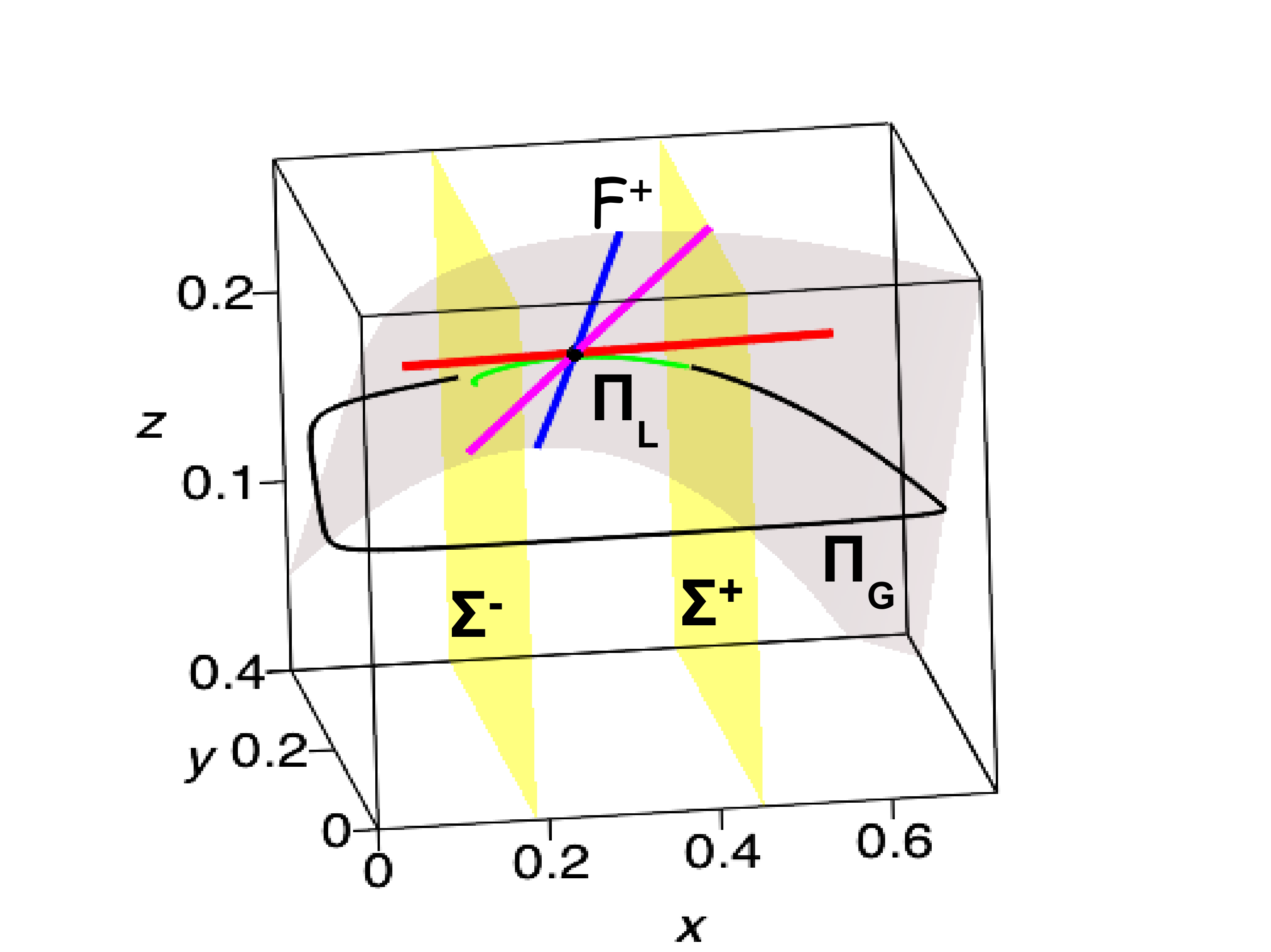} }}%
		\end{center}
		 \caption{(a): The critical manifold $\mathcal{M}=\mathcal{S}\cap\mathcal{T}$.
     (b): The folded node is represented 
		by a black dot, the strong and weak eigendirections are represented 
		by red and magenta lines respectively. The singular funnel is the 
		region between the blue and the red curves.     \footnote{Reprinted from 
		\cite{SCT1}, Copyright (2017) with permission from Elsevier.}
		}
 \label{global_return}%
\end{figure}

With the assumption $\zeta\ll 1$, we obtain a singularly
perturbed system in one fast and two slow variables which can be written
on the fast time scale as
\begin{eqnarray}
\label{nondim2}    
\left\{
\begin{array}{ll} x'&= x\left(1-x-\frac{y}{\beta_1+x}
-\frac{z}{\beta_2+x}\right):=xu(x,y,z),\\
   y'&=\zeta y\left(\frac{x}{\beta_1+x}-c\right)
	:=\zeta yv(x,y,z),\\
     z' &= \zeta  z\left(\frac{x}{\beta_2+x}-d-hz\right)
		:=\zeta zw(x,y,z),
       \end{array} 
\right. \end{eqnarray}
where $u(x,y,z)$, $v(x,y,z)$, $w(x,y,z)$ represent the 
nontrivial $x$, $y$ and $z$ nullclines.  Hence,
we may apply techniques from fast-slow systems~\cite{KuehnBook}.  
The {\em{critical manifold}} 
$\mathcal{M}$ is the set of equilibria of the fast subsystem, i.e.
\bess 
\mathcal{M}=\left\{(x,y,z): x=0 {\mathrm{~or~}} u(x,y,z)=0\right\} 
=: \mathcal{T}\cup\mathcal{S}, 
\eess
where $\mathcal{T} =\{ (0,y,z): y, z \geq 0\}$, 
$\mathcal{S}=\{(x,y,z):u(x,y,z)=0\}$. The critical manifold
consists of two normally attracting sheets $\mathcal{S}^a$ and $\mathcal{T}^a$, and two 
repelling sheets $\mathcal{S}^r$ and $\mathcal{T}^r$, separated by two curves of 
singularities $\mathcal{F}^+$ and $\mathcal{F}^-$ \cite{S2, SCT1}.
The normally hyperbolic parts perturb to slow manifolds by 
Fenichel's Theorem~\cite{Fenichel4,KuehnBook}. The 
curve $\mathcal{F}^-$ of transcritical points is the intersection of 
$\mathcal{T}$ with $\mathcal{S}$. The curve $\mathcal{F}^+=\{(x,y,z): u(x,y,z)=0
,\ u_x(x,y,z)=0\}$ is a curve of fold points. On rescaling $s:=\zeta t$, 
\eqref{nondim2} transforms to an equivalent system  written on the
slow time scale $s$. On letting $\zeta \to 0$ on this scale, one obtains 
the slow subsystem (or reduced subsystem). The corresponding dynamics 
is referred to as the slow flow which is restricted to the 
plane $\{x=0\}$ or to the surface $\mathcal{S}:=\{(x,y,z):u(x,y,z)=0\}$. 
Since $u_y \neq 0$ on $\mathcal{S}$, the implicit function theorem implies that 
$S$ can be locally written as a graph of $y=\phi(x,z)$. Differentiating 
$u(x,y,z)=0$ implicitly with respect to slow time yields
\bess
u_x\dot{x}+u_y\dot{y}+u_z\dot{z}=0,\qquad 
\frac{\textnormal{d}}{\textnormal{d}s}=\dot{~}. 
\eess
Consequently, the slow flow restricted to $\mathcal{S}$ 
reads as
\begin{eqnarray} 
\label{red2}  
\begin{pmatrix}
 -u_x\dot{x}  \\
 \dot{z} 
 \end{pmatrix}= \begin{pmatrix}
 u_yyv+u_zzw  \\
zw
 \end{pmatrix}\bigg|_{y=\phi(x,z)}.
\end{eqnarray}
Multiplying the vector field by a factor $-u_x$, and then using a suitable 
time re-scaling, transforms~\eqref{red2} to the desingularized system
\begin{eqnarray} 
\label{dsing}  \begin{pmatrix}
 \dot{x}  \\
 \dot{z} 
 \end{pmatrix}= \begin{pmatrix}
 u_yyv+u_zzw  \\
-u_xzw
 \end{pmatrix}\bigg|_{y=\phi(x,z)}.
\end{eqnarray}
The set of all points on $\mathcal{F}^+$ for which $u_yyv+u_zzw \neq 0$ 
corresponds to jump points and are said to satisfy the normal 
switching condition~\cite{Benoit2,SW,KuehnBook}. At these points, solutions 
jump from $\mathcal{S}^a$ to $\mathcal{T}^a$. Similarly, one can treat jump points at the curves 
of transcritical singularities in $\mathcal{T}$~\cite{KruSzm4,KuehnSzmolyan}. This 
gives rise to relaxation-oscillation dynamics as shown in 
Fig.~\ref{global_return}(b). Points on $\mathcal{F}^+$ that violate 
the normal switching condition are known as folded singularities~\cite{Benoit2,
SW,KuehnBook}. These are equilibria of the desingularized 
system~\eqref{dsing}. Folded singularities can be classified into their
equilibrium type, e.g., one generically obtains folded nodes, folded 
foci and folded saddles without parameter variation; further degeneracies
occur, e.g., folded-saddle nodes~\cite{KrupaWechselberger,W1} and/or singular 
Hopf bifurcation~\cite{Braaksma,Guckenheimer7}. 

Under suitable parameter values, the desingularized system~\eqref{dsing} 
admits a folded node singularity $(x^{\star},y^{\star},z^{\star})$. Suppose 
that $\lambda_s, \lambda_w$ are the eigenvalues of the Jacobian 
corresponding to~\eqref{dsing} at $(x^{\star},y^{\star},z^{\star})$ 
such that $\lambda_s< \lambda_w<0$. The perturbation for $0<\zeta\ll1$ 
of trajectories tangent to the eigendirections for $\zeta=0$ 
corresponding to eigenvalues $\lambda_s, \lambda_w$ gives rise to 
primary canards, known as the strong canard, and the weak canard 
denoted by $\gamma_s$ and 
 $\gamma_w$, respectively~\cite{Benoit2,SW,W1,KuehnBook}. The singular funnel 
refers to the region on $\mathcal{S}^{a}$ bounded by the strong canard and 
the fold curve $\mathcal{F}^+$. It contains trajectories flowing 
through the folded node. For $0<\zeta\ll 1$, these trajectories 
experience small rotations around the weak canard until they jump 
to the other attracting sheet of the manifold $\mathcal{T}^a$. A global return 
mechanism can re-inject trajectories to the folded-node funnel to 
induce MMOs. The global return mechanism 
can be represented by a return map $\Pi: \Sigma^+ \to \Sigma^+$, 
where $\Sigma^+$ is a cross-section orthogonal to the $x$-axis away 
from the fold $\mathcal{F}^+$ such that for small $\zeta$, all 
trajectories that are projected into the slow manifold $\mathcal{S}^a_{\zeta}$
onto a neighborhood of the funnel intersect with $\Sigma^+$. 
As shown in Fig.~\ref{global_return}(b), the map $\Pi$ can be
 decomposed as 
$\Pi = \Pi_G\circ\Pi_L,$ 
where $\Pi_L:\Sigma^+\to \Sigma^-$ defines the local map near the 
fold (shown in green) and $\Pi_G: \Sigma^-\to \Sigma^+$ the global return 
map (in black)~\cite{DGKKOW,K}. We will study a stochastic 
variant of one part of this return map in Section~\ref{sec:main}.   
  
Suppose we are in a parameter regime such that~\eqref{nondim2} 
has a folded node singularity at 
$(x^{\star},y^{\star}, z^{\star})$. Consider the eigenvalue ratio 
\begin{equation}
\label{eq:evratio}
\mu:=\lambda_w/\lambda_s.
\end{equation}
Then one may conclude~\cite{SW, W1,KuehnBook}, for sufficiently small 
$\zeta>0$, the strong singular canard $\gamma_{s}$ perturbs 
to a maximal canard $\gamma_s^{\zeta}$ and if $\mu^{-1}\notin 
\mathbb{N}$, the singular weak canard $\gamma_w$ perturbs to a 
maximal canard $\gamma_w^{\zeta}$. Moreover if $k\in \mathbb{N}$ 
be such that 
\bess
2k-1<\mu^{-1}<2k+1\quad \textnormal{ and }\quad \mu^{-1}\neq 2(k+1),
\eess
there exists $k$ other additional canards referred to as secondary 
canards. In this case, the strong singular canard $\gamma_s$ twists 
around the primary weak canard $\gamma_w^{\zeta}$. The $j$-th secondary 
canard $\gamma_j^{\zeta}$, $1\leq j\leq k$ twists $2j+1$ times around 
$\gamma_w^{\zeta}$. The twisting occurs in an $O(\sqrt{\zeta})$-neighborhood 
of the folded node  $(x^{\star},y^{\star}, z^{\star})$. This mechanism 
induces small-amplitude oscillations (SAOs).  A related mechanism of SAO
generation is the case, when a folded saddle-node of type 
II occurs in system (\ref{dsing}), or more commonly, a singular Hopf bifurcation occurs 
in system (\ref{nondim2}) ~\cite{SW,Guckenheimer7,KuehnBook}.  
In this case, SAOs are induced by dynamics near an equilibrium point of 
the full system with complex conjugate eigenvalues~\cite{DGKKOW}. To illustrate
this in our context, we fix the parameter values to 
\bes
\label{par_values}
\textnormal{$\zeta=0.01, \ \beta_1=0.5,\  \beta_2=0.25,\ c=0.38,\ d=0.17$}
\ees
and view $h$ as the main bifurcation parameter. Then as discussed in~\cite{S2},
 system~\eqref{nondim2} 
exhibits a Hopf bifurcation at 
\begin{equation}
\label{eq:Hopf}
h\approx 2.6413\qquad \textnormal{(Hopf of the full system)} 
\end{equation}
and the desingularized 
system~\eqref{dsing} undergoes a folded saddle-node bifurcation of type II 
at 
\begin{equation}
\label{eq:FSNII}
h\approx 2.722\qquad \textnormal{(FSN II in singular limit)}.
\end{equation} 
Hence, the deterministic model shows many
classical features, which have been studied in other fast-slow system
models, mainly in chemical reactions and neuroscience~\cite{KuehnBook}. 


\section{The Stochastic Model}
\label{sec:main}

In this section, we will build and analyze a stochastic variant 
of~\eqref{nondim2}. Our first step is to justify the structure of 
the noise terms from a microscopic description while modeling the birth 
and death processes of the species as continuous time Markov processes 
on a discrete state space. We will consider large population sizes so 
that we can apply the law of large numbers allowing us to use continuous 
variables and differential equations in a mean-field limit~\cite{Gardiner}. 
A similar approach has been taken in~\cite{CL} for modeling the interactions 
between a predator and its prey, where the authors consider a hybrid 
model treating the predator as a continuous variable and the 
prey as discrete. In our work, the evolution of the three species 
are assumed to take discrete values. To this end, let $\omega_1 x(t)$ 
be an integer denoting the number of prey at time $t$. In other words, 
$x$ units of prey corresponds to $\omega_1 x$ individuals. Similarly, 
let $\omega_2 y(t)$ and $\omega_3 z(t)$ be integers denoting the number 
of two different predator species at time $t$.  

Suppose that $\omega_1x$ follows an exponential distribution with 
parameter $\lambda_x$ with expectation $1/\lambda_x$, where 
 $\lambda_x ={\omega_1}/{\zeta} .$
Then the number of events $N_x$ in terms of births and deaths 
during $[s, s+\textnormal{d}s]$ is approximately 
\bess
N_x\approx \frac{\textnormal{d}s}{(1/\lambda_x)}=
\frac{\omega_1}{\zeta} \textnormal{d}s.
\eess
Let $X_i$ be the random variable representing a change in the prey 
population at the $i$-th event. Suppose that $X_i$ equals $1$ if 
a birth occurs, $-1$ if a death occurs and $0$ if no changes occur  
with the following probabilities:
\begin{eqnarray*}
\mathbb{P}(X_i=1) &=&\frac{x}{1+x+x^2+\frac{xy}{\beta_1+x}+
\frac{xz}{\beta_2+x}},\\
\mathbb{P}(X_i=-1)&=&\frac{x^2+\frac{xy}{\beta_1+x}+
\frac{xz}{\beta_2+x}} {1+x+x^2+\frac{xy}{\beta_1+x}+\frac{xz}{\beta_2+x}},\\
 \mathbb{P}(X_i=0)&=& \frac{1}{1+x+x^2+\frac{xy}{\beta_1+x}+
\frac{xz}{\beta_2+x}}.
\end{eqnarray*}
Then the expectation and the variance of the random variables 
$X_i$ and $\sum_{i=1}^{N_x} X_i$ read as
\begin{eqnarray*}
 \mathbb{E}(X_i) &=& \frac{x - x^2- \frac{xy}{\beta_1+x}- 
\frac{xz}{\beta_2+x}}{1+x+x^2+\frac{xy}{\beta_1+x}+\frac{xz}{\beta_2+x}}\\
 \mathbb{E}(\sum_{i=1}^{N_x} X_i) &=& \frac{\omega_1}{\zeta}
\Big( x - x^2- \frac{xy}{\beta_1+x}- \frac{xz}{\beta_2+x}\Big)\textnormal{d}s, \\
 \sigma^2(X_i) &=& \mathbb{E}(X_i^2)- (\mathbb{E}(X_i))^2, \\
 \sigma^2(\sum_{i=1}^{N_x} X_i) &=& \frac{\omega_1}{\zeta}F_1^2(x,y,z) 
~\textnormal{d}s,
\end{eqnarray*}
where 
\begin{eqnarray}
 \label{Fn1} F_1^2(x,y,z) &= &\frac{  x + \Big(x^2+ \frac{xy}{\beta_1+x}
+ \frac{xz}{\beta_2+x}\Big)(1+4x)}{1+x+x^2+\frac{xy}{\beta_1+x}
+\frac{xz}{\beta_2+x}} .
\end{eqnarray}
Assuming that $\omega_1$ is sufficiently large, from the central 
limit theorem, $\sum_{i=1}^{N_x} X_i$ can be approximated by a 
Gaussian random variable:
\[ \sum_{i=1}^{N_x} X_i \approx \frac{\omega_1}{\zeta}f_1(x,y,z) 
\textnormal{d}s + \sqrt{\frac{\omega_1}{\zeta}} F_1(x,y,z) \  W_s^{x}, \]
where $W_t^{x}$ is the standard one-dimensional Wiener process, 
$F_1(x,y,z)$ is the positive square root of~\eqref{Fn1} and
\begin{eqnarray}
\label{fn1} f_1(x,y,z) &=& x - x^2- \frac{xy}{\beta_1+x}- \frac{xz}{\beta_2+x}.
\end{eqnarray}
 Since the 
variable $x$ represents the number of individuals of prey divided 
by $\omega_1$, the increment in $x$ during $[s, s+\textnormal{d}s]$ will 
be thus given by
\begin{eqnarray}
 \label{incx} x(s+\textnormal{d}s)-x(s) =  
\frac{f_1(x,y,z)}{\zeta}\textnormal{d}t + 
\frac{F_1(x,y,z)}{\sqrt{\omega_1 \zeta}}  \  W_s^{x}.
\end{eqnarray}
Similarly, assume that $\omega_2 y$ and $\omega_3z$ follow 
exponential distributions with parameters $\lambda_y$ and 
$\lambda_z$ respectively, where $\lambda_y= \omega_2$ and 
$\lambda_z=\omega_3$. Let $Y_i$ and $Z_i$ be the random variables 
representing changes in the populations of the predators at the 
$i$-th event taking values $\pm 1$ or $0$ with probabilities given
in Table~\ref{table:prob}.
 
\begin{table}[h]
\centering 
\begin{tabular}{c  c} 
\hline\hline 
Random variable & Probability  \\ [0.5ex] 
\hline 
$ Y_i=1$  & $\frac{\frac{xy}{\beta_1+x}}{1+
\frac{xy}{\beta_1+x}+cy}$ \\ 
 $Y_i=-1$ & $ \frac{cy} {1+\frac{xy}{\beta_1+x}+cy} $ \\
$Y_i=0$ & $ \frac{1} {1+\frac{xy}{\beta_1+x}+cy}$  \\
$Z_i=1$  & $\frac{\frac{xz}{\beta_2+x}}{1+\frac{xz}{\beta_2+x} 
+dz+ hz^2}$   \\
$Z_i=-1$& $ \frac{dz+hz^2}{1+\frac{xz}{\beta_2+x} +dz+ hz^2}$  \\ 
$Z_i=0$ & $ \frac{1}{1+\frac{xz}{\beta_2+x} +dz+ hz^2}$  
\\ [1ex] 
\hline 
\end{tabular}
\caption{Probabilities associated with the changes} 
\label{table:prob} 
\end{table}

Then similar to the derivation of increment of $x$, the increments of 
$y$ and $z$ in the time interval $[s, s+\textnormal{d}s]$ are given 
by the equations
\begin{eqnarray}
\label{incy} y(s+\textnormal{d}s)-y(s) & = & f_2(x,y,z)\textnormal{d}s 
+ \frac{F_2(x,y,z)}{\sqrt{\omega_2}}  \  W_s^{y},\\
 \label{incz} z(s+\textnormal{d}s)-z(s)  &= & f_3(x,y,z) \textnormal{d}s 
+ \frac{F_3(x,y,z)}{\sqrt{\omega_3}} \  W_s^{z},
\end{eqnarray}
 respectively, where 
\begin{eqnarray}
\label{fn2} f_2(x,y,z) &=&   \frac{xy}{\beta_1+x}- cy, \\
\label{fn3}f_3(x,y,z) &=& \frac{xz}{\beta_1+x}- dz-hz^2, \\
\label{Fn2} F_2(x,y,z) &=& \sqrt{\frac{ \frac{xy}{\beta_1+x}
+cy +\frac{4cxy^2}{\beta_1+x}}{1+\frac{xy}{\beta_1+x}+cy}}, \\
\label{Fn3} F_3(x,y,z) &=&   \sqrt{\frac{\frac{xz}{\beta_2+x}
+\Big(1+\frac{4xz}{\beta_2+x}\Big)(dz+hz^2)}{1
+\frac{xz}{\beta_2+x}+dz+hz^2}}.
\end{eqnarray}
Denoting $\sigma_i = 1/\sqrt{\omega_i}$ for $i=1,2,3$ and 
$X_t$, $Y_t$ and $Z_t$ as the random variables for the three 
interacting populations, combining equations~\eqref{incx}-\eqref{incz}, 
we thus obtain a system of SDEs interpreted in 
It$\hat{\mathrm{o}}$-sense:
\begin{eqnarray}\label{stocorig}    \left\{
\begin{array}{ll}  \textnormal{d}X_s =&  \frac{1}{\zeta} 
f_1(X_s, Y_s, Z_s) \ \textnormal{d}s+\frac{\sigma_1}{\sqrt{\zeta}} 
F_1(X_s, Y_s, Z_s) \ \textnormal{d}W_s^{(1)}\\
\textnormal{d}Y_s =& f_2(X_s, Y_s, Z_s)\ \textnormal{d}s+\sigma_2 
F_2(X_s, Y_s, Z_s) \ \textnormal{d}W_s^{(2)}\\
\textnormal{d}Z_s =& f_3(X_s, Y_s, Z_s)\ \textnormal{d}t+\sigma_3 
F_3(X_s, Y_s, Z_s)\ \textnormal{d}W_s^{(3)},
 \end{array} 
\right. \end{eqnarray}
where $f_i, F_i$ are given by equations~\eqref{Fn1}-\eqref{fn1} 
and~\eqref{fn2}-\eqref{Fn3} for $i=1,2,3$. 

\begin{figure}[ht!]%
    \centering
       \subfloat[$h=2.66$]{{\includegraphics[width=4.0cm]
			{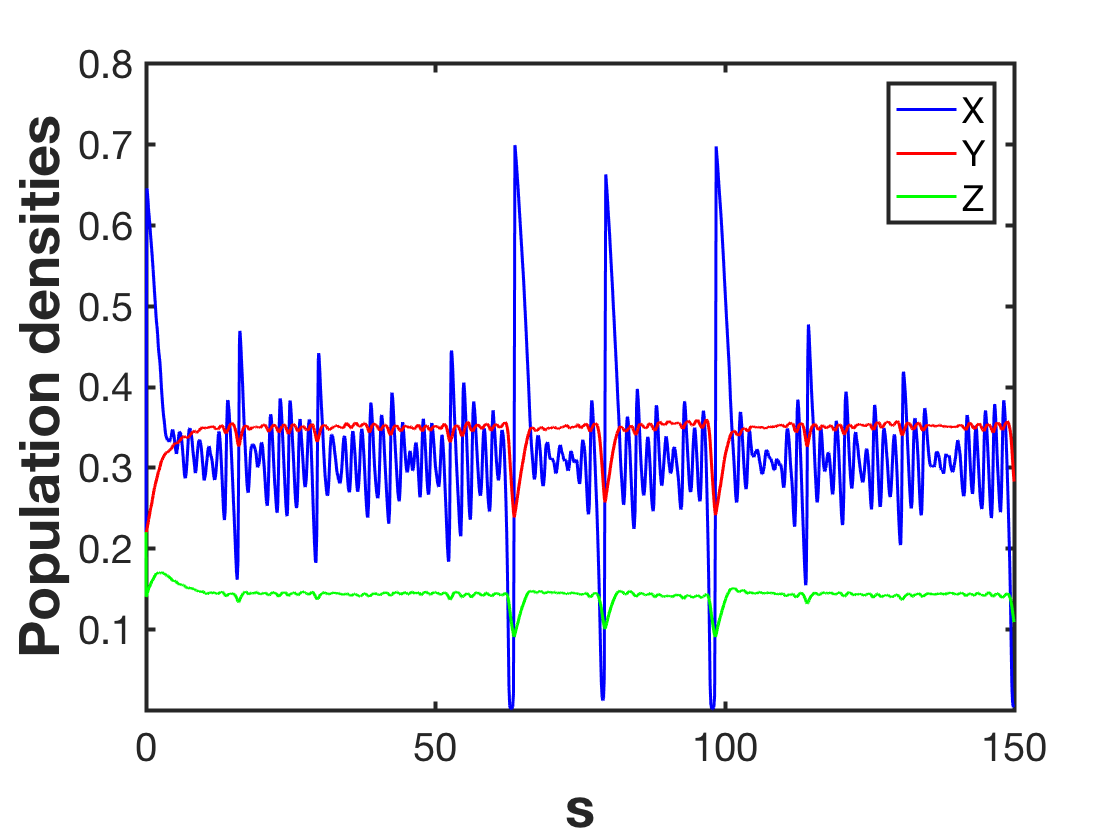} }}%
           \quad
           \subfloat[$h=2.66$]{{\includegraphics[width=4.0cm]
					{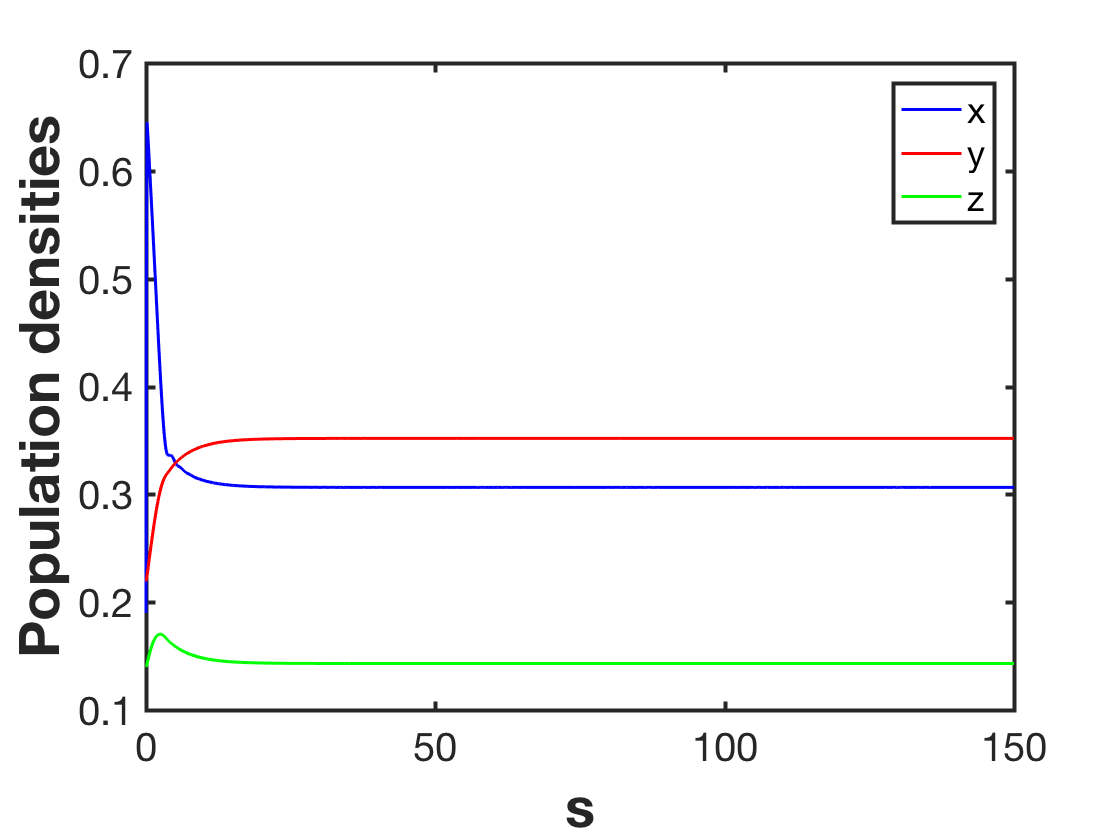} }}%
               \quad
                   \subfloat[$h=2.6$]{{\includegraphics[width=4.0cm]
									{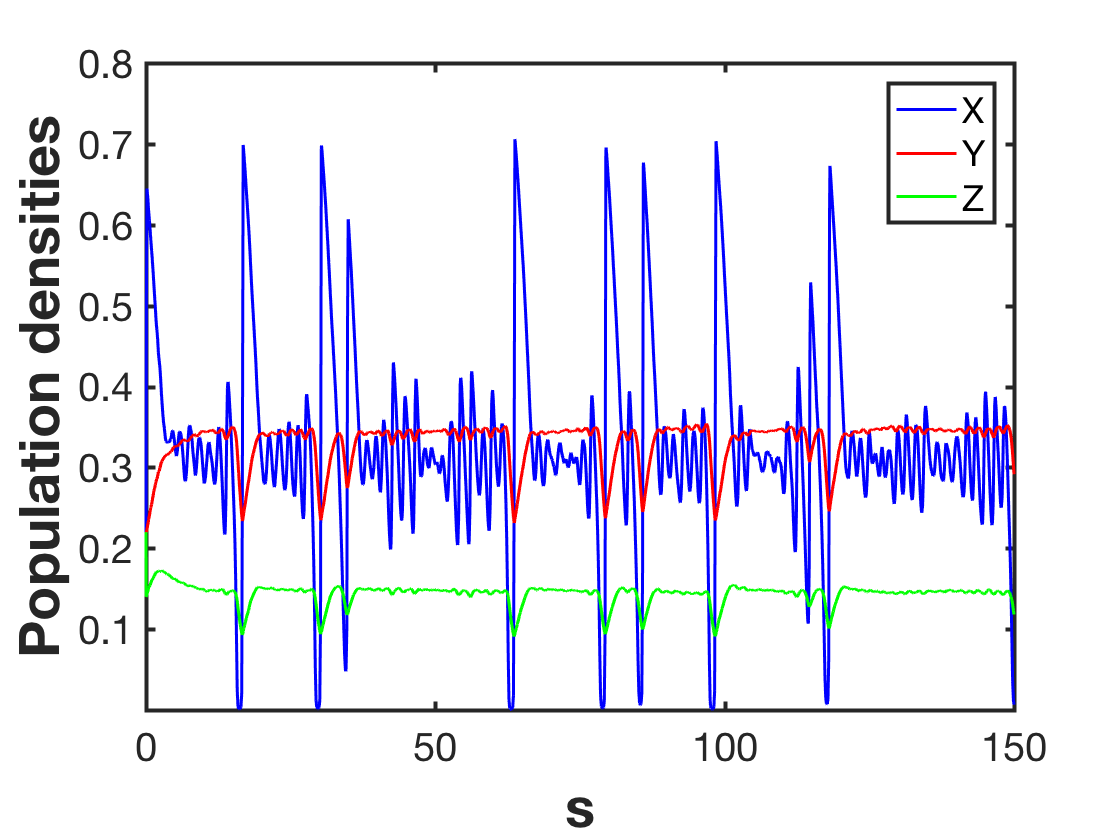} }}%
                       \quad
     \subfloat[$h=2.6$]{{\includegraphics[width=4.0cm]
		{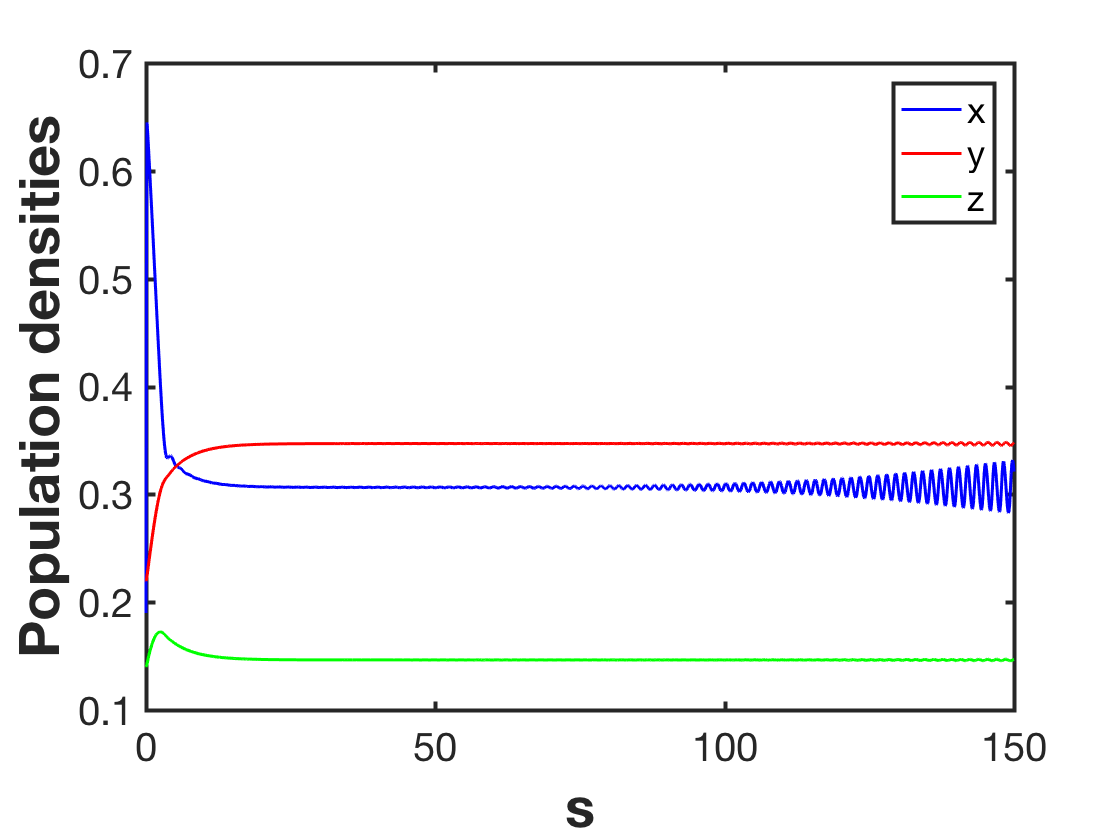} }}%
  \quad
    \subfloat[$h=2.4$]{{\includegraphics[width=4.0cm]
		{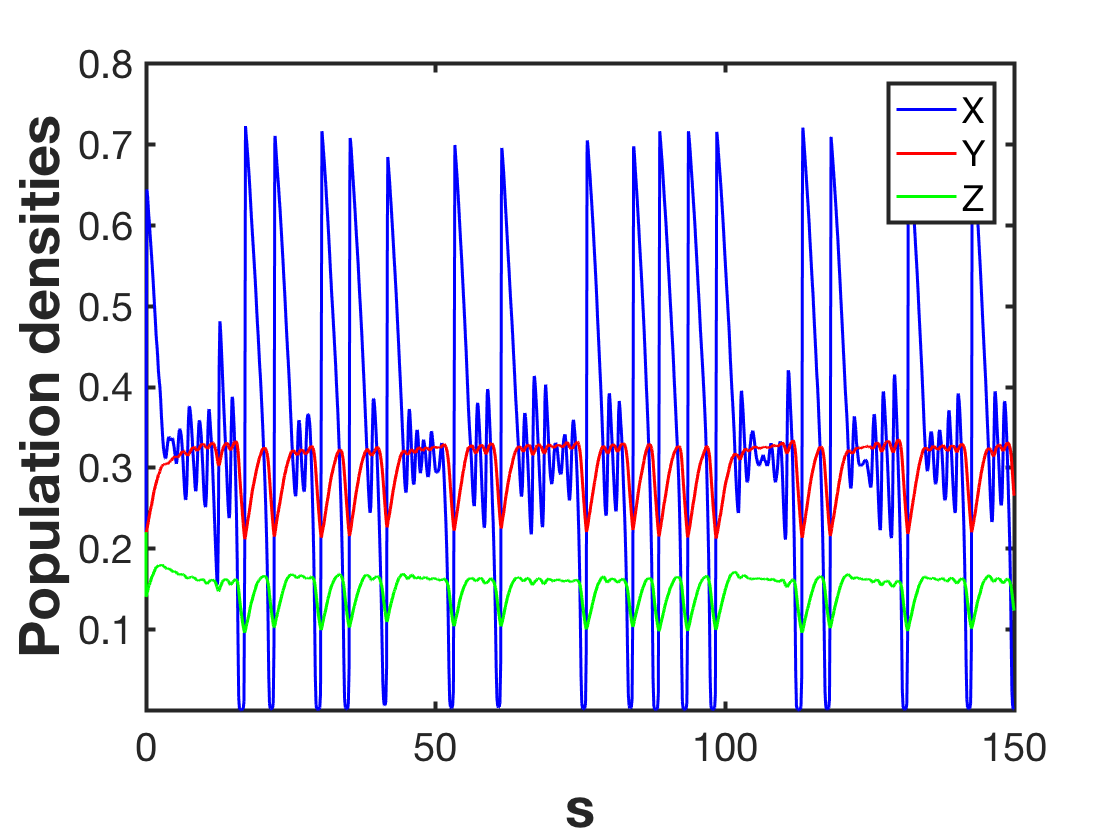} }}%
    \quad
    \subfloat[$h=2.4$]{{\includegraphics[width=4.0cm]
		{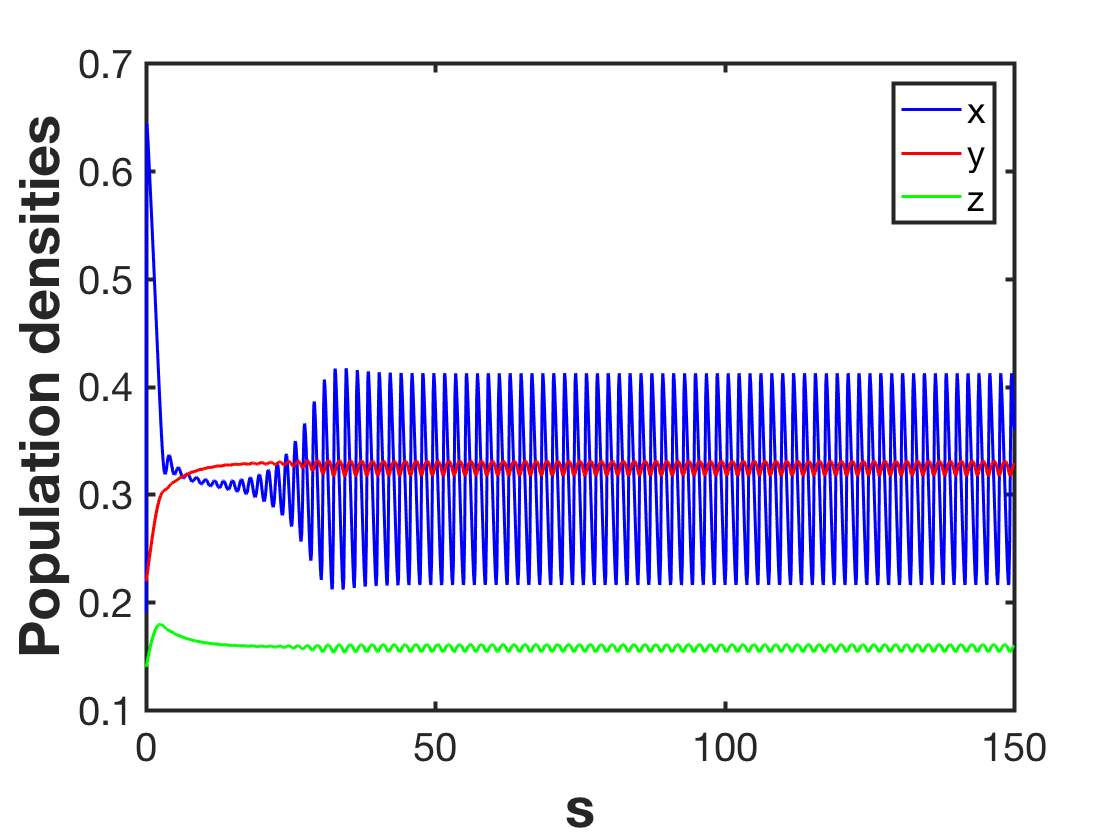} }}%
\caption{Left column: Time series of different sample paths 
of system~\eqref{stocorig} as $h$ is varied. Note the existence 
of isolated spikes and repeated spikes.
Right column: Time series of the deterministic 
system~\eqref{nondim2} on the slow time scale for varying $h$.
The other parameter values are given by (\ref{par_values})
with $\sigma_1 = 10^{-6},$ 
$\sigma_2=\sigma_3=3\cdot10^{-3}$.}
    \label{timeseries_0}
\end{figure}

Rescaling $s= \zeta t$ and using the scaling law of Brownian motion, 
the equivalent fast time scale system of SDEs reads as
 \begin{eqnarray}\label{stoc}    \left\{
\begin{array}{ll} \textnormal{d}X_t =&  f_1(X_t, Y_t, Z_t) 
\textnormal{d}t+\sigma_1F_1(X_t, Y_t, Z_t) \textnormal{d}W_t^{(1)}\\
\textnormal{d}Y_t =& \zeta f_2(X_t, Y_t, Z_t)\textnormal{d}t
+\sqrt{\zeta}\sigma_2F_2(X_t, Y_t, Z_t) \textnormal{d}W_t^{(2)}\\
 \textnormal{d}Z_t =& \zeta f_3(X_t, Y_t, Z_t) \textnormal{d}t
+\sqrt{\zeta}\sigma_3  F_3(X_t, Y_t, Z_t) \textnormal{d}W_t^{(3)}.
 \end{array} 
\right. \end{eqnarray}
Systems~\eqref{stocorig} and~\eqref{stoc} are referred to as 
stochastic fast-slow systems with $X$ being the fast and $Y$, $Z$ 
being the slow variables, and $\sigma_i$ for $i=1,2,3$ are the 
corresponding noise parameters. 

Note that $\sigma_i^2$ is inversely proportional to the size of 
the population. For larger population size we can expect the 
deterministic part of the process to dominate on quite long time 
scales~\cite{BerglundGentz}. Consequently, on these time scales, 
extremely small values of $\sigma_i$ (for instance, $\sigma_1 = 
10^{-7}, \sigma_2=\sigma_3=10^{-4}$) lead to only very small 
deviations from the deterministic model~\eqref{nondim2}.

\begin{figure}[ht!]%
    \centering
       \subfloat[$h=2.69$.]{{\includegraphics[width=4.0cm]
			{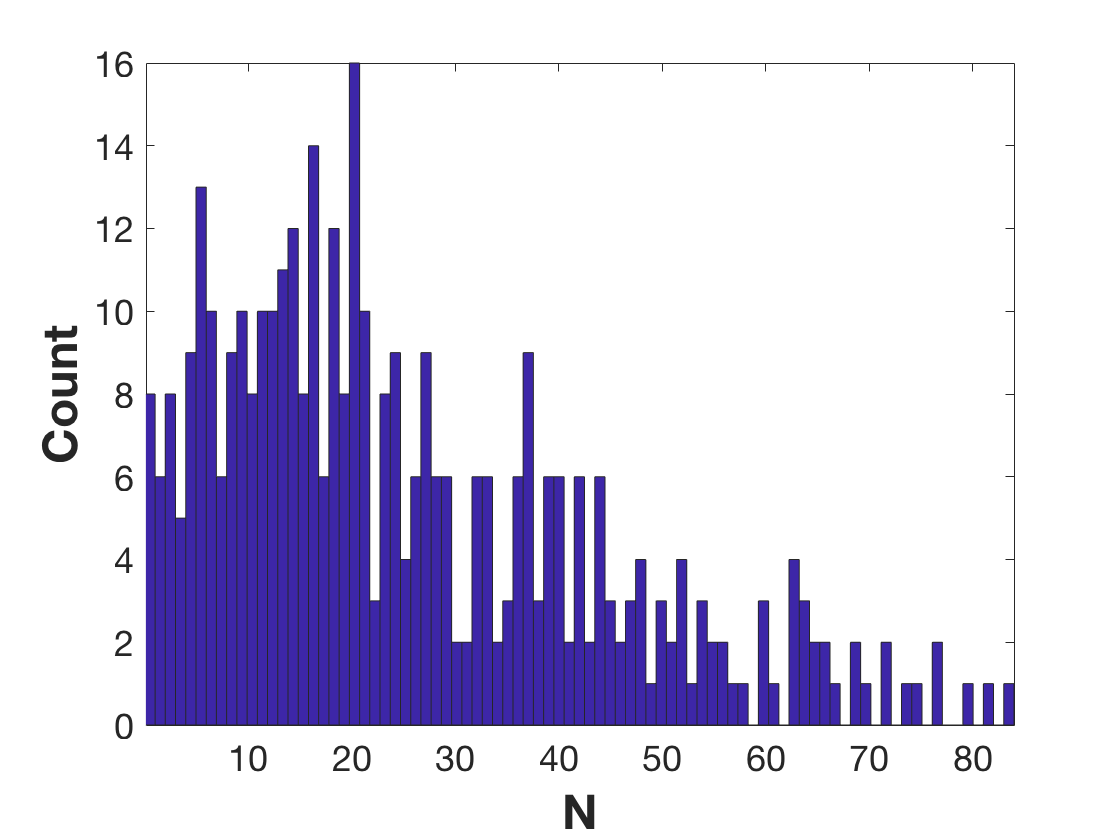} }}%
   \quad
    \subfloat[$h=2.66$.]{{\includegraphics[width=4.0cm]
		{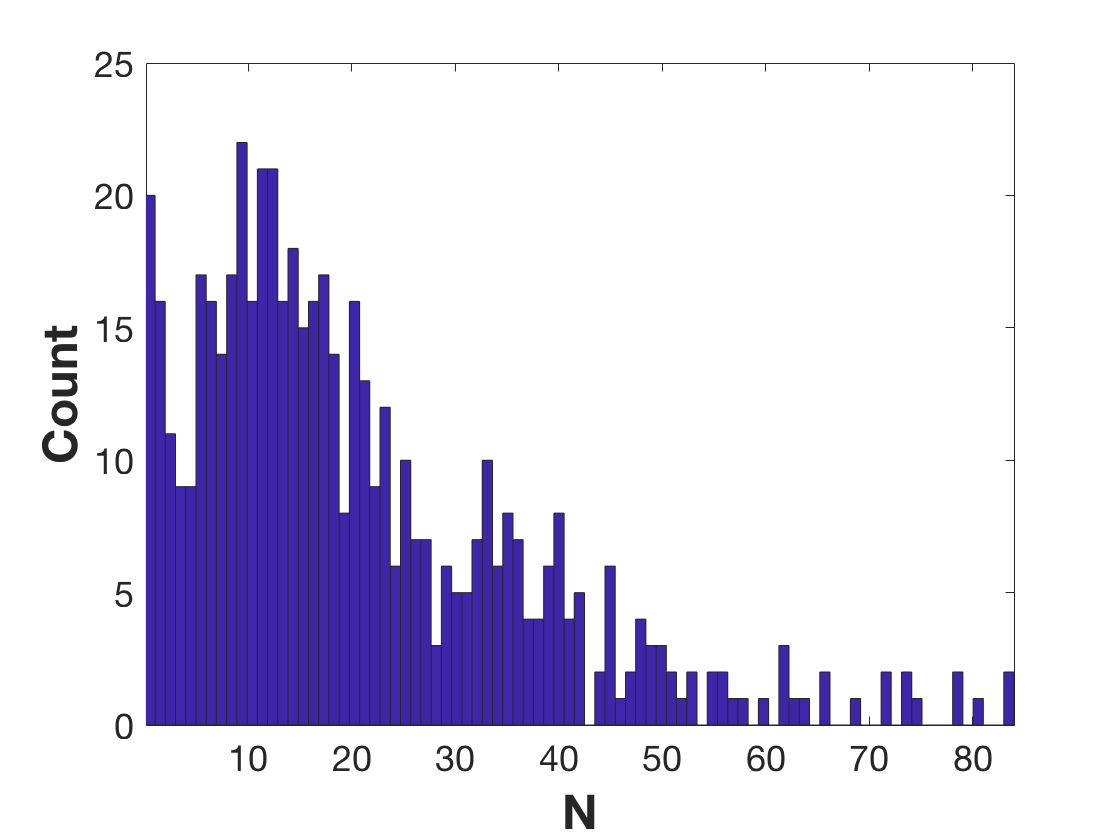}} }%
      \quad
    \subfloat[$h=2.64$.]{{\includegraphics[width=4.0cm]
		{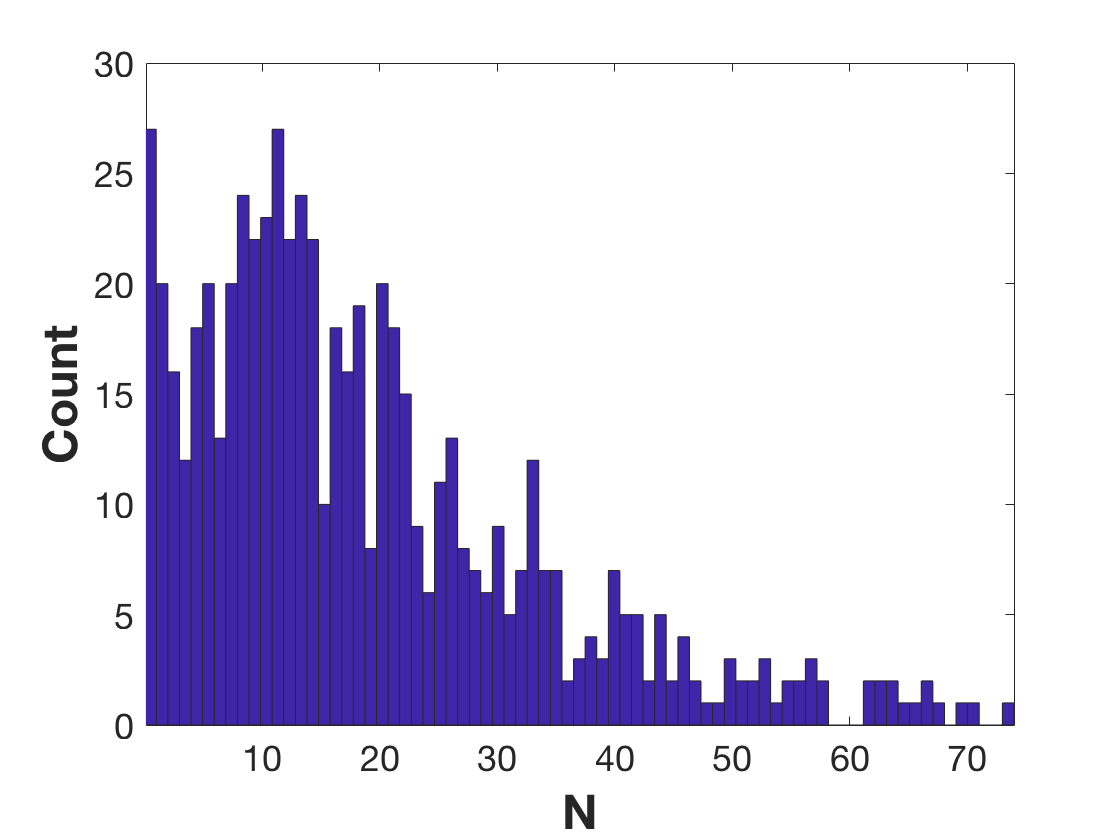} }}%
    \quad
    \subfloat[$h=2.6$.]{{\includegraphics[width=4.0cm]
		{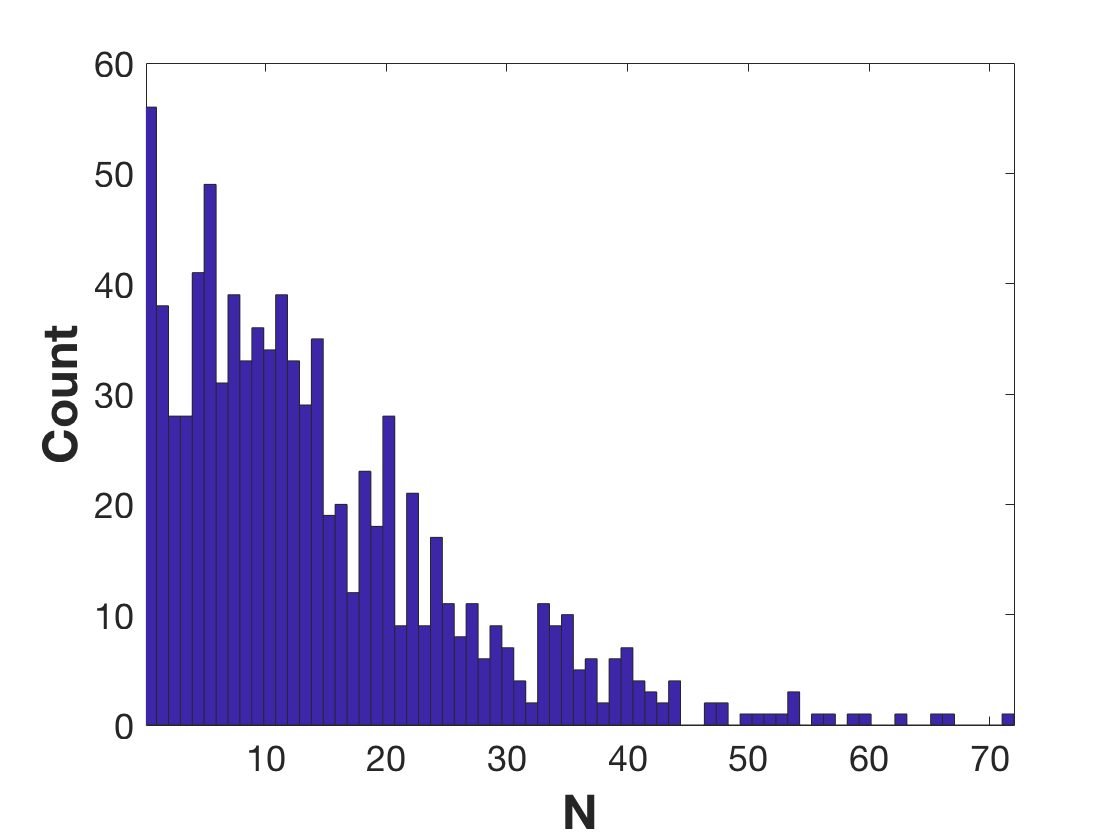}}}
  \quad
        \subfloat[$h=2.5$.]{{\includegraphics[width=4.0cm]
				{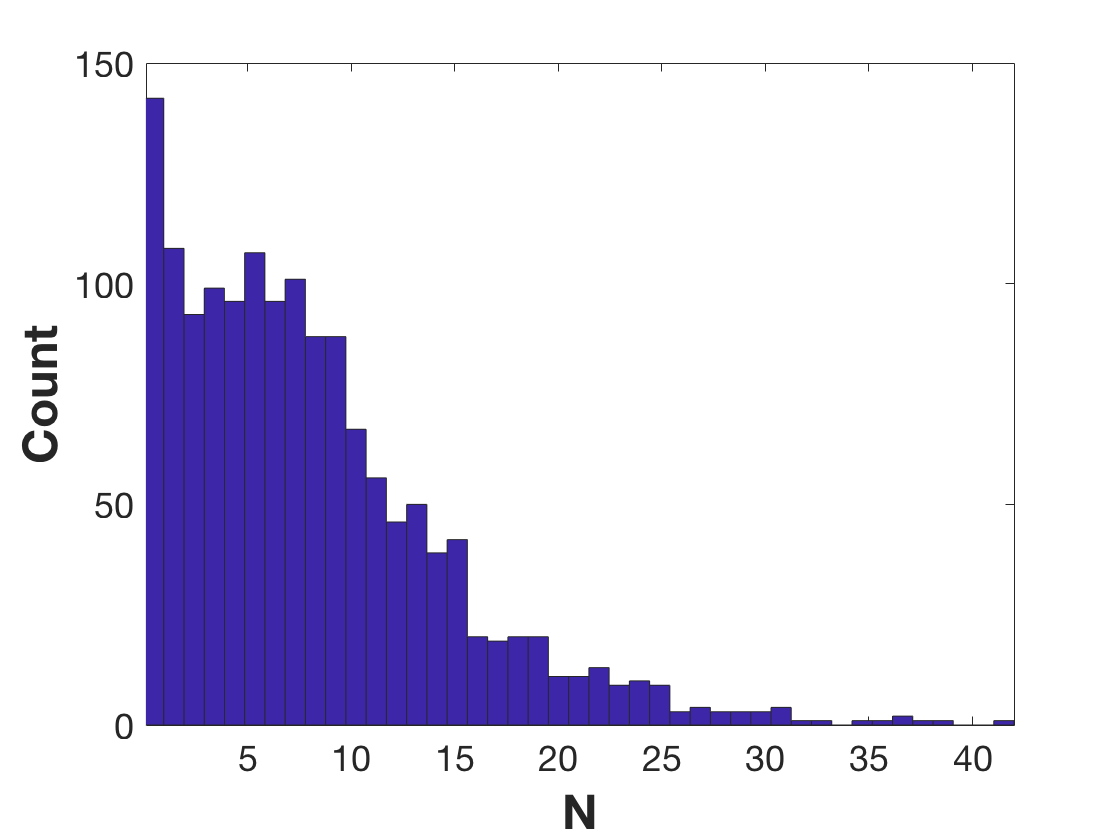} }}%
    \quad
    \subfloat[$h=2.4$.]{{\includegraphics[width=4.0cm]
		{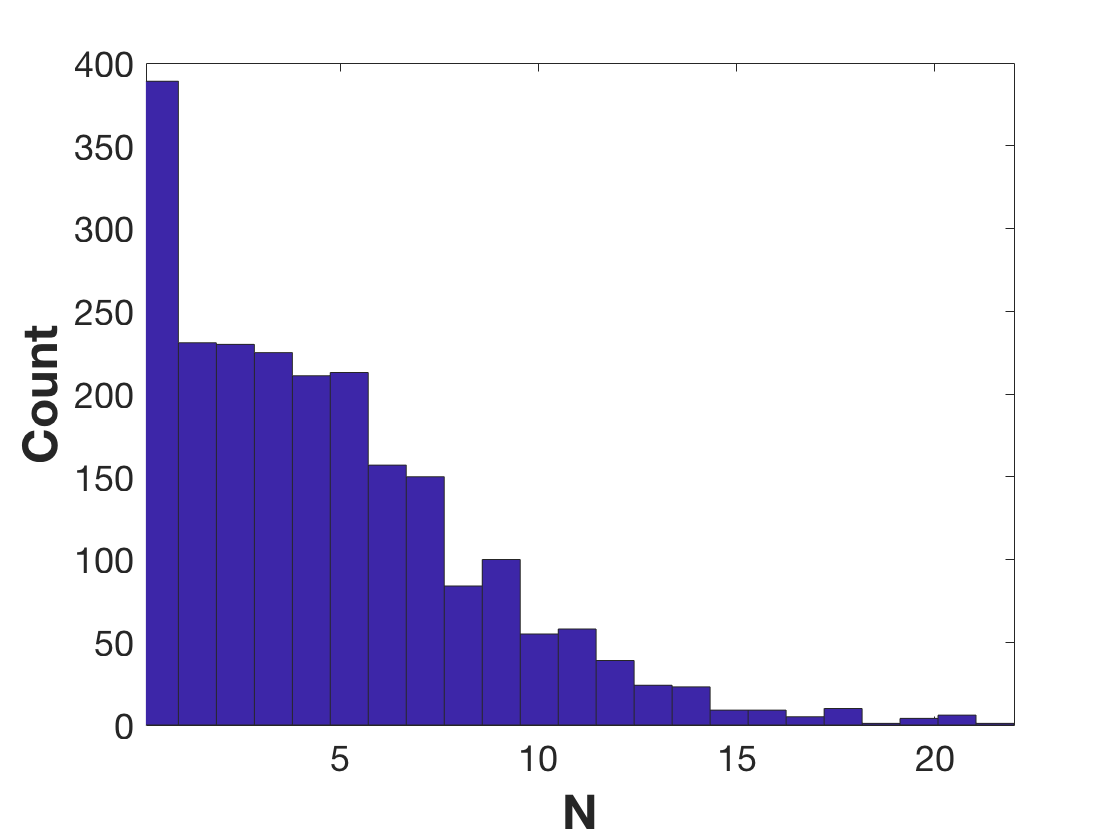}}}
\caption{Histograms of the distribution of the number of SAOs 
 between two spikes corresponding to
system~\eqref{stocorig} close to Hopf 
bifurcation for parameter values given by (\ref{par_values})
 and noise 
intensities $\sigma_1= 10^{-6}$, $\sigma_2 = \sigma_3= 3\cdot10^{-3}$. 
    \label{histogram_original_1}}%
\end{figure}

\begin{figure}[ht!]%
    \centering
      \subfloat[$\sigma_1=10^{-6},\sigma_2=10^{-3},
			\sigma_3=10^{-5}$.]{{\includegraphics[width=4.0cm]
			{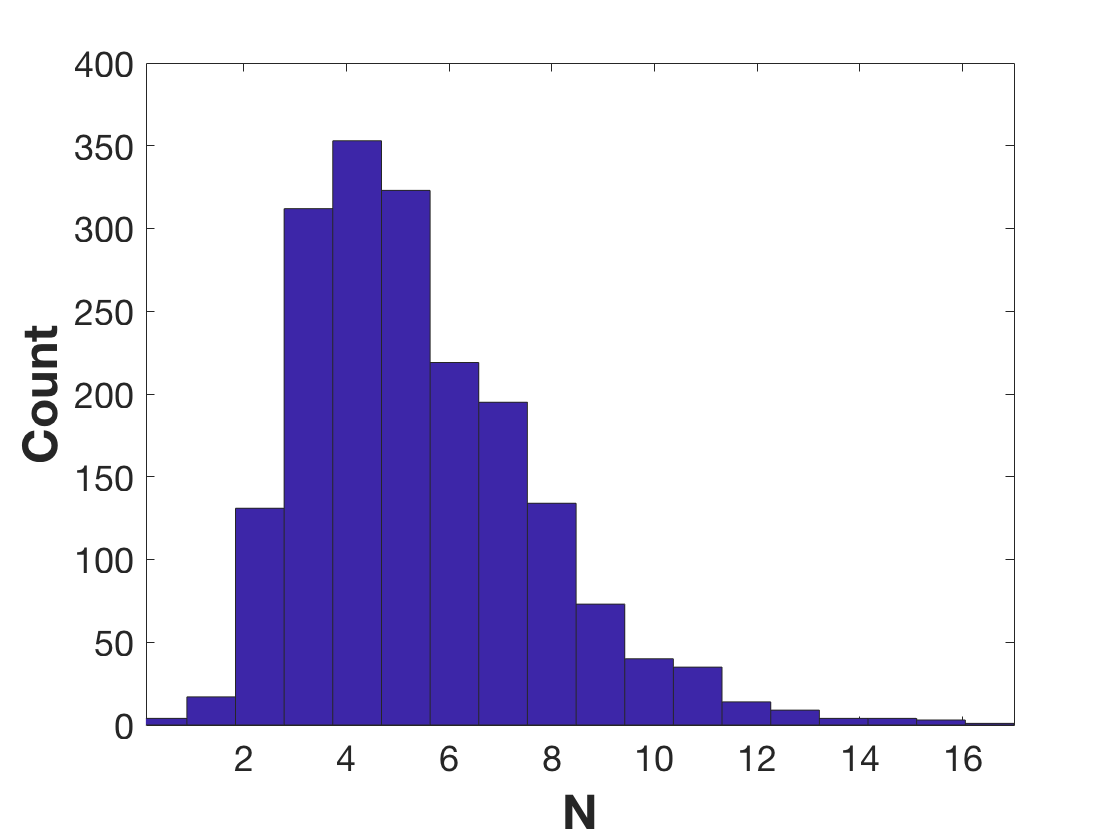} }}%
      \quad
    \subfloat[$\sigma_1=10^{-6},\sigma_2=10^{-3},
		\sigma_3=10^{-4}$]{{\includegraphics[width=4.0cm]
		{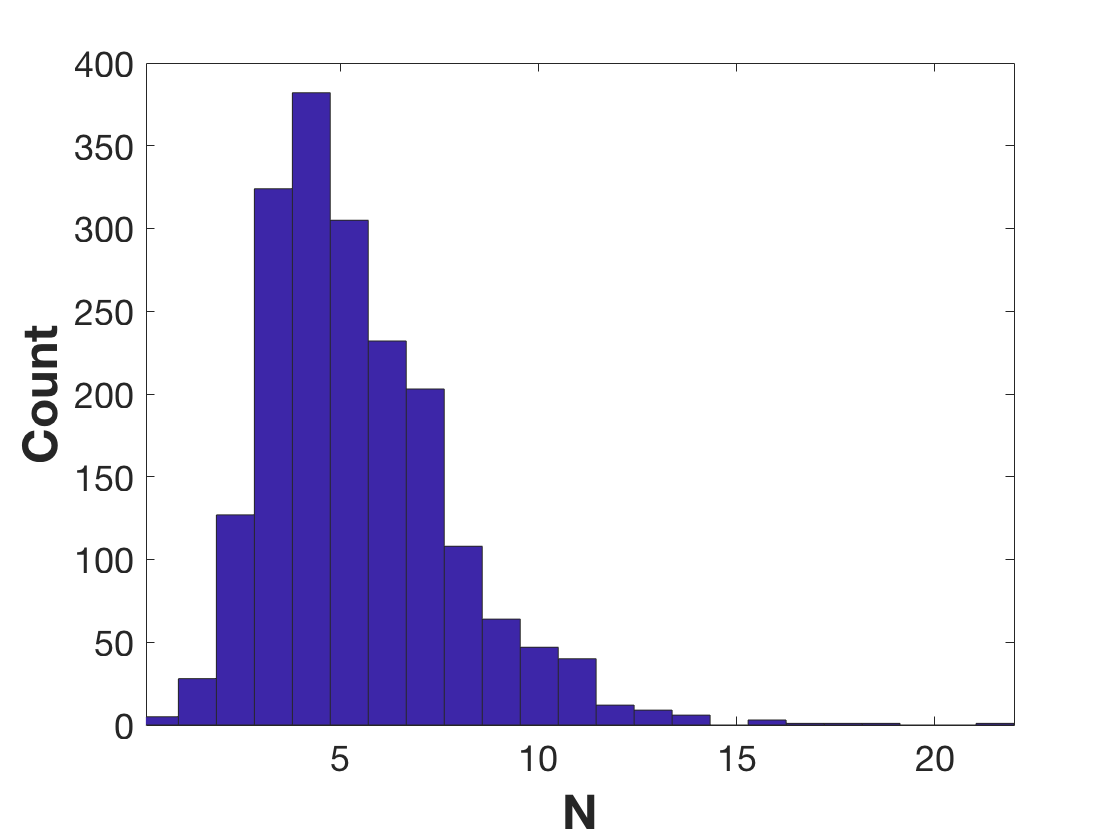} }}%
    \quad
   \subfloat[$\sigma_1=10^{-6},\sigma_2=10^{-3},
	\sigma_3=10^{-3}$]{{\includegraphics[width=4.0cm]
	{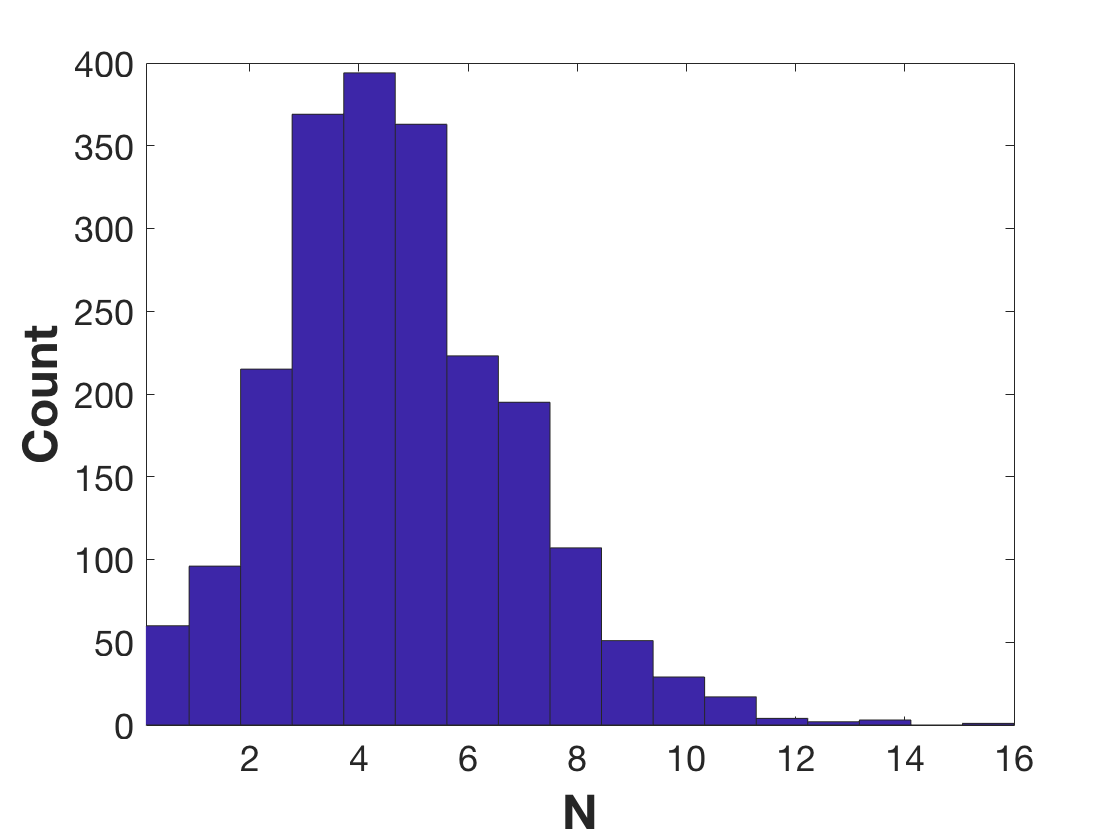} }}
    \quad
    \subfloat[$\sigma_1=10^{-6},\sigma_2=10^{-3},
		\sigma_3=5\cdot10^{-3}$]{{\includegraphics[width=4.0cm]
		{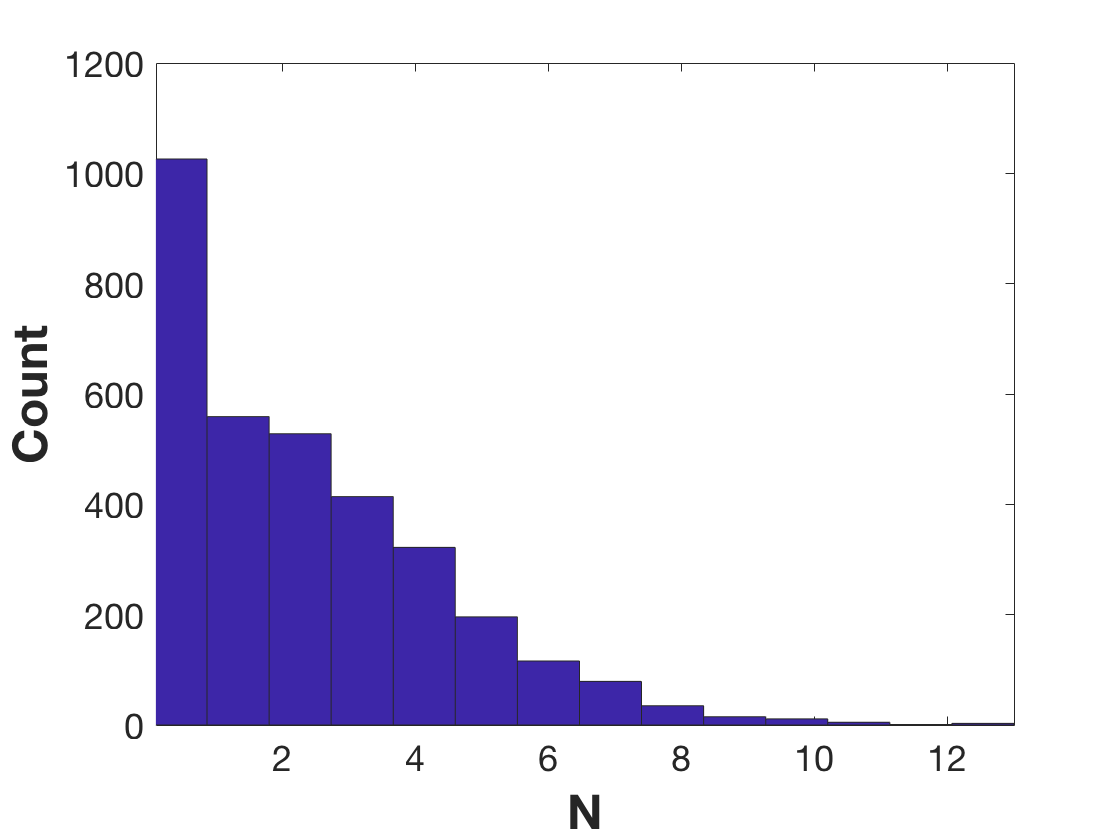} }}
\caption{Histograms of SAOs corresponding to system~\eqref{stocorig} 
at $h=2.3$ for different noise intensities and parameter values given by (\ref{par_values}).}
    \label{hist_stochastic_1}%
\end{figure}

Larger values of $\sigma_i$ (such as $\sigma_1=10^{-5}$, 
$\sigma_2=\sigma_3=10^{-4}$) can lead to extinction of species. This is 
easy to see from our derivation as the noise terms vanish when the population 
sizes are zero. 

Interestingly, for intermediate values of $\sigma_i$  (for 
instance, $\sigma_1 = 10^{-6}, \sigma_2=\sigma_3=3\cdot10^{-4}$), we observe different kinds 
of noise-induced MMOs as shown in the numerical
simulations~\cite{Higham} in Fig.~\ref{timeseries_0}. The isolated/repeated spikes 
 in the time series in Fig.~\ref{timeseries_0} represent 
events of random population outbreaks or collapses as documented for small 
mammals and insects in the wild~\cite{AC, M1, KBJP}. We also observe 
that the noise-driven dynamics can be significantly 
different depending on how close we are to the singular Hopf bifurcation;
 cf.~\eqref{eq:Hopf}-\eqref{eq:FSNII}. Indeed, if we are below it, we   
see rare excursions giving rise to LAOs (isolated spikes). On the other hand, in the near Hopf-regime, 
due to the small stable limit cycle generated 
at the (singular) Hopf bifurcation, which has a part near an unstable 
slow manifold~\cite{DGKKOW}, the noise can induce LAOs more frequently. 
This means we really see noise-generated MMOs with two different characteristics.

\begin{figure}[ht!]%
    \centering
			\includegraphics[width=7.0cm]{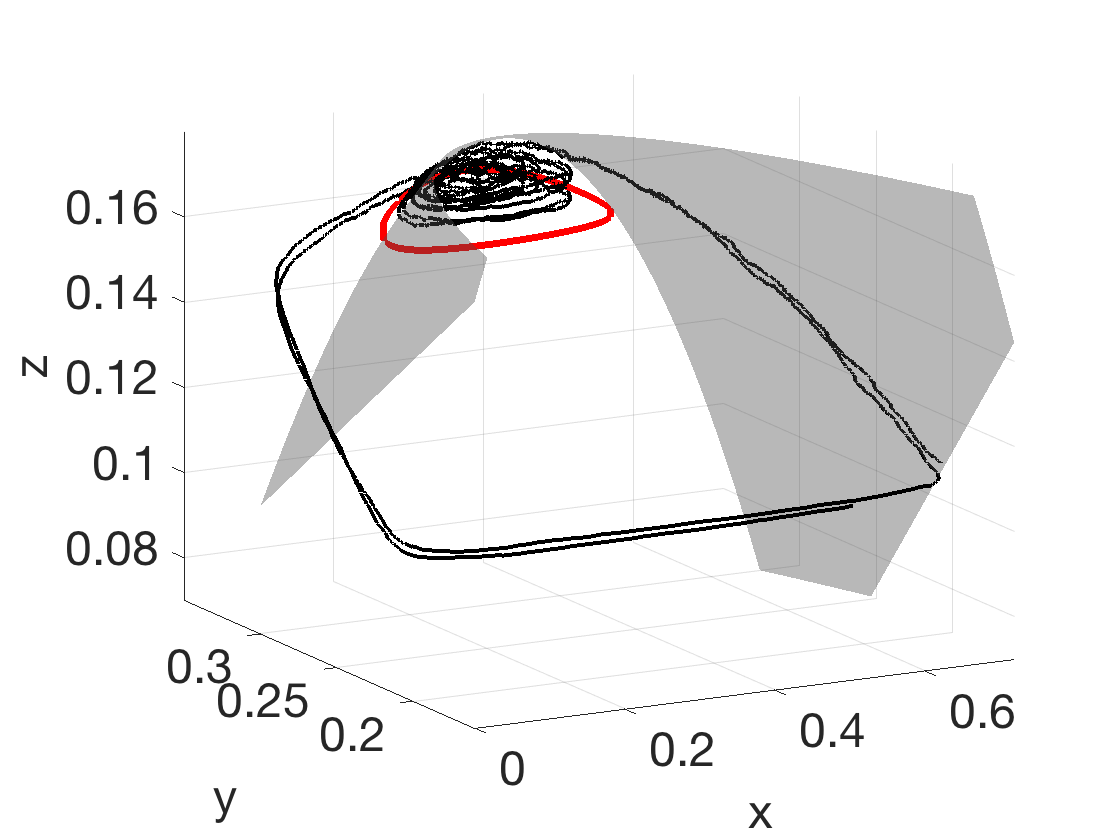}%
  \caption{A sample path of system (\ref{stocorig}) along with 
  the deterministic stable limit cycle (plotted in red) against the critical manifold
  $\mathcal{S}$.
  The other parameter values are as in Fig.~\ref{hist_stochastic_1} with 
  $\sigma_1=10^{-6},$ $\sigma_2=\sigma_3=5\cdot10^{-3}$.} 
  \label{3d_samplepath}%
\end{figure}

 Fig.~\ref{timeseries_0} also suggests the presence of
  an intricate interplay between the noise-induced SAOs and the SAOs initiated 
by the drift terms. To investigate this effect further, we consider the
number of SAOs between two spikes and denote it by $N$. The distribution of $N$ for $200$ 
different sample paths near the Hopf bifurcation are 
shown in Fig.~\ref{histogram_original_1} and Fig.~\ref{hist_stochastic_1}; 
in our numerical study, we considered oscillations in the $x$-variable.
Oscillations whose amplitudes exceeded $0.1$ with 
their maximum value below $0.6$ in the $x$-coordinate were considered 
as SAOs. We observe that for parameter values above and close to the 
singular Hopf bifurcation, the distributions are multimodal with wider range,
 while  below and away from the Hopf value, the distributions are reverse J-shaped 
 with shorter range. Similarly, lower noise levels lead to
 unimodal distributions with long tails,  while higher noise yield distributions that  are
reverse J-shaped (see Fig.~\ref{hist_stochastic_1}). We conjecture that the reverse J-shaped distributions
  asymptotically approach an exponential distribution with increasing noise levels or
  distance from the Hopf bifurcation. 
Similar results have been shown in a planar model with multiplicative noise in \cite{S1}. 
In fact, the noise-induced SAO count depends on the number of oscillations
a sample path makes around 
the deterministically-stable small limit cycle.  A large excursion occurs
only if the sample path gets kicked by noise from a neighborhood of 
this limit cycle as shown in Fig.~\ref{3d_samplepath}.
We refer to~\cite{BerglundLandon,MuratovVanden-Eijnden} for the simpler case of 
additive noise in planar systems, which is already highly non-trivial
mathematically. We also remark that the histograms in Fig.~\ref{histogram_original_1} 
are naturally related to ecological applications, such as the 
distribution of return times of outbreaks of larch budmoth in Euopean 
Alps~\cite{EBFNL}. To make this connection
more precise would be an interesting direction for future work.

\begin{figure}[ht!]%
    \centering
       \subfloat[$(\sigma_1,\sigma_2,\sigma_3)=(10^{-6},10^{-3},
		10^{-4})$]{{\includegraphics[width=4.0cm]
			{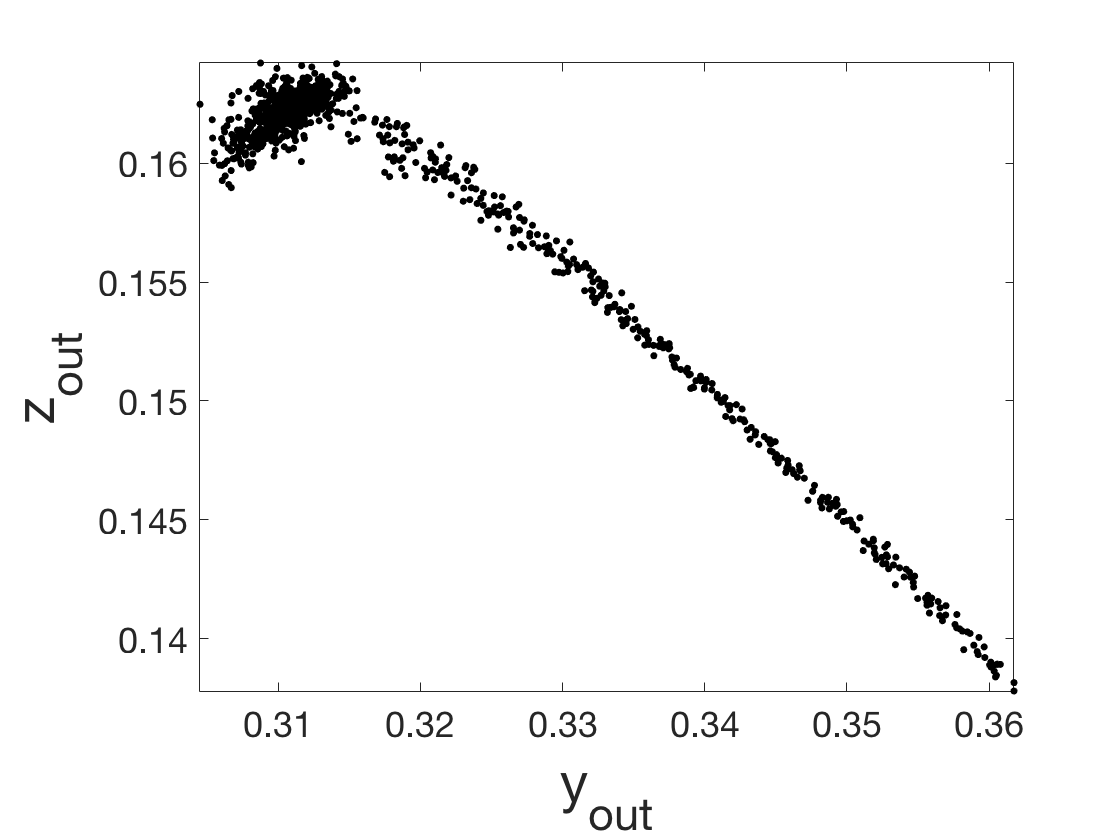} }}%
   \quad
    \subfloat[$(\sigma_1,\sigma_2,\sigma_3)=(10^{-6},10^{-3},
		10^{-4})$]{{\includegraphics[width=4.0cm]
		{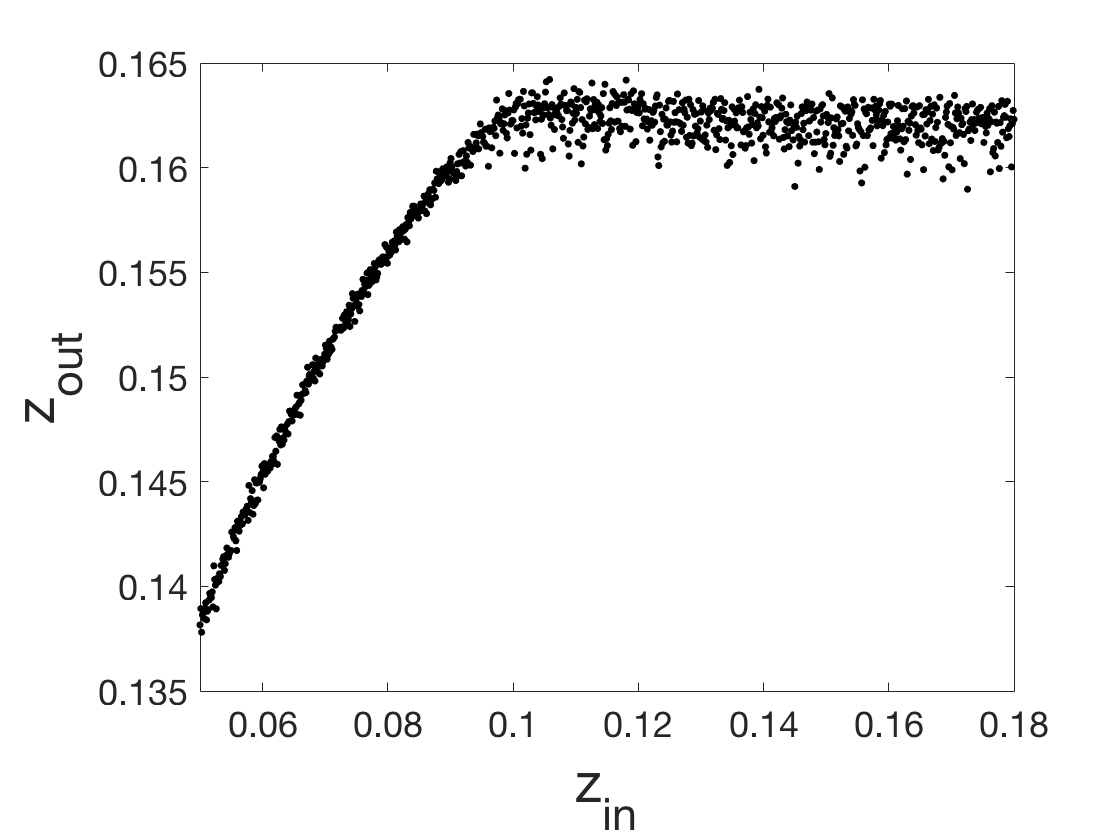} }}%
\quad    
    \subfloat[$(\sigma_1,\sigma_2,\sigma_3)=(10^{-6},10^{-3},
		10^{-3})$]{{\includegraphics[width=4.0cm]
		{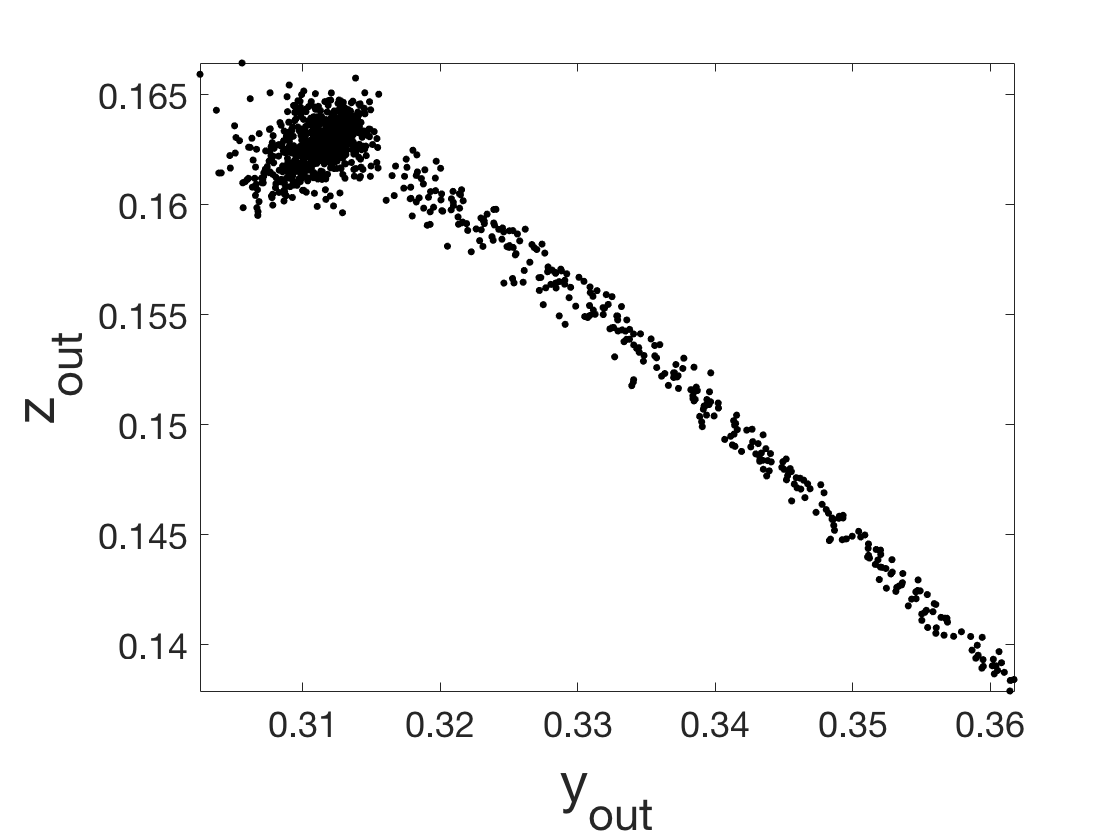} }}%
    \quad
    \subfloat[$(\sigma_1,\sigma_2,\sigma_3)=(10^{-6},10^{-3},
		10^{-3})$]{{\includegraphics[width=4.0cm]
		{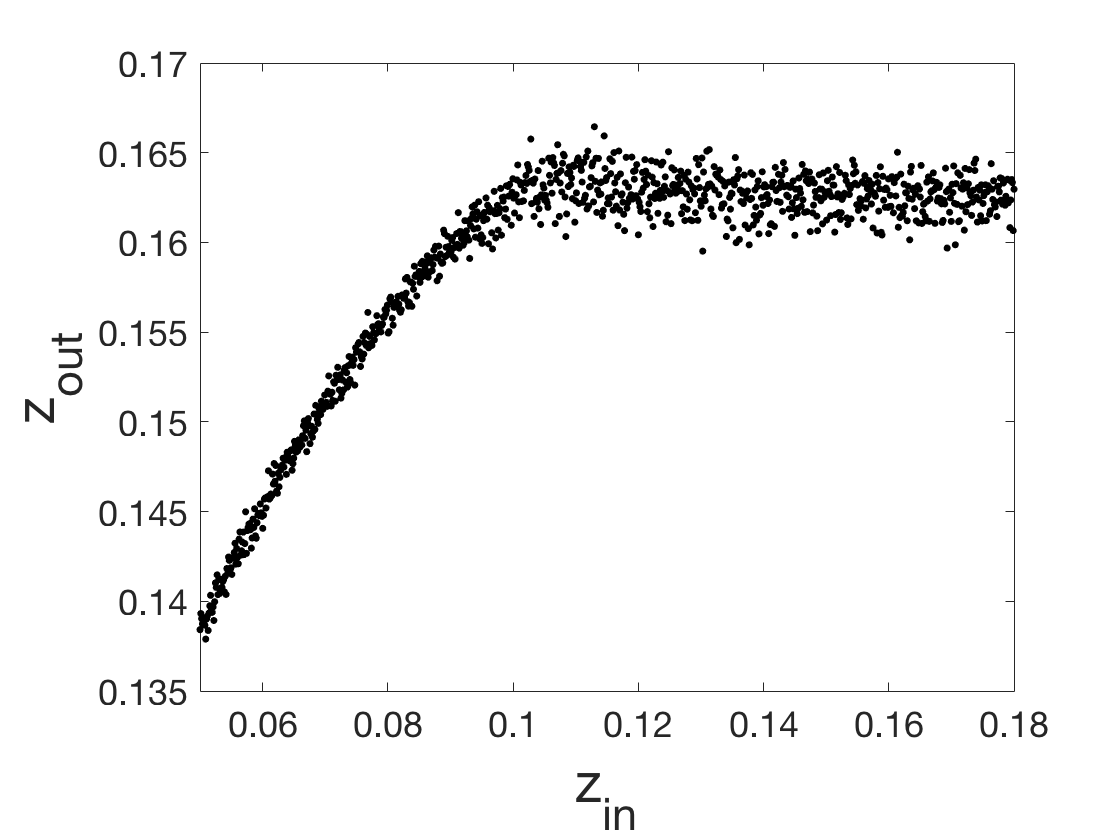} }}%
\quad
     \subfloat[$(\sigma_1,\sigma_2,\sigma_3)=(10^{-6},10^{-3},
		10^{-2})$.]{{\includegraphics[width=4.0cm]
		{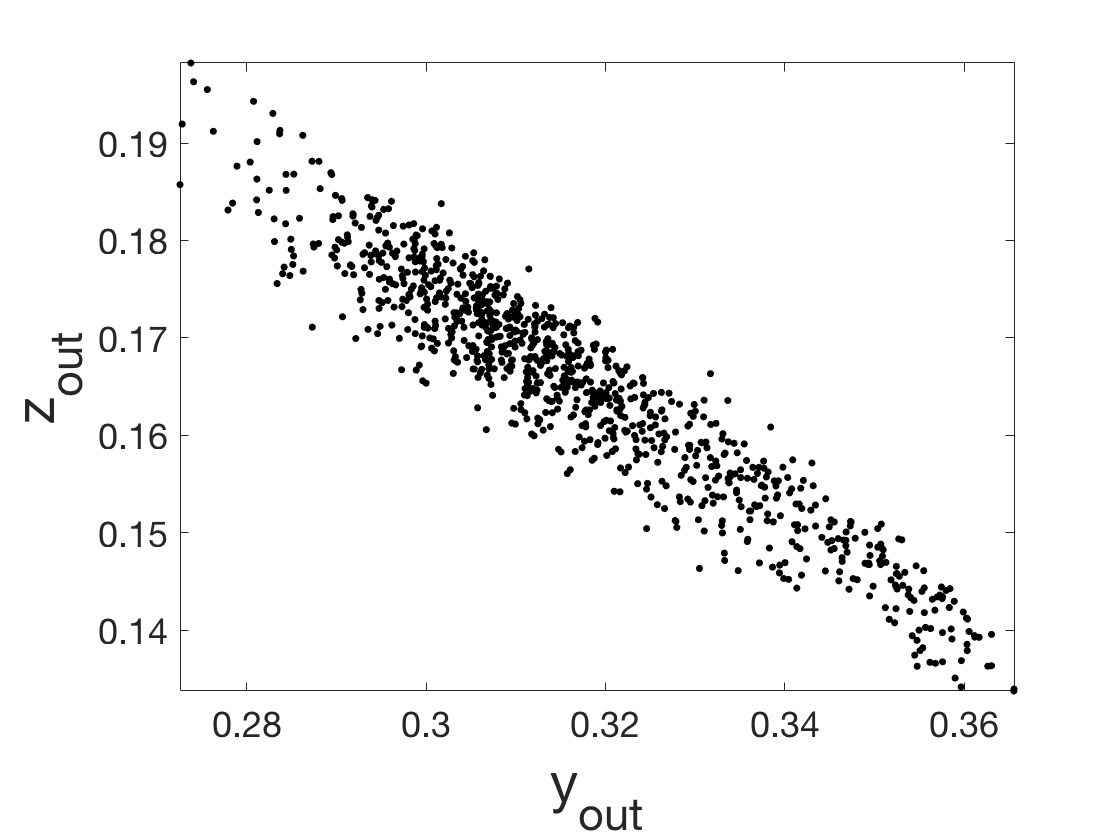} }}%
    \quad
    \subfloat[$(\sigma_1,\sigma_2,\sigma_3)=(10^{-6},10^{-3},
		10^{-2})$]{{\includegraphics[width=4.0cm]
		{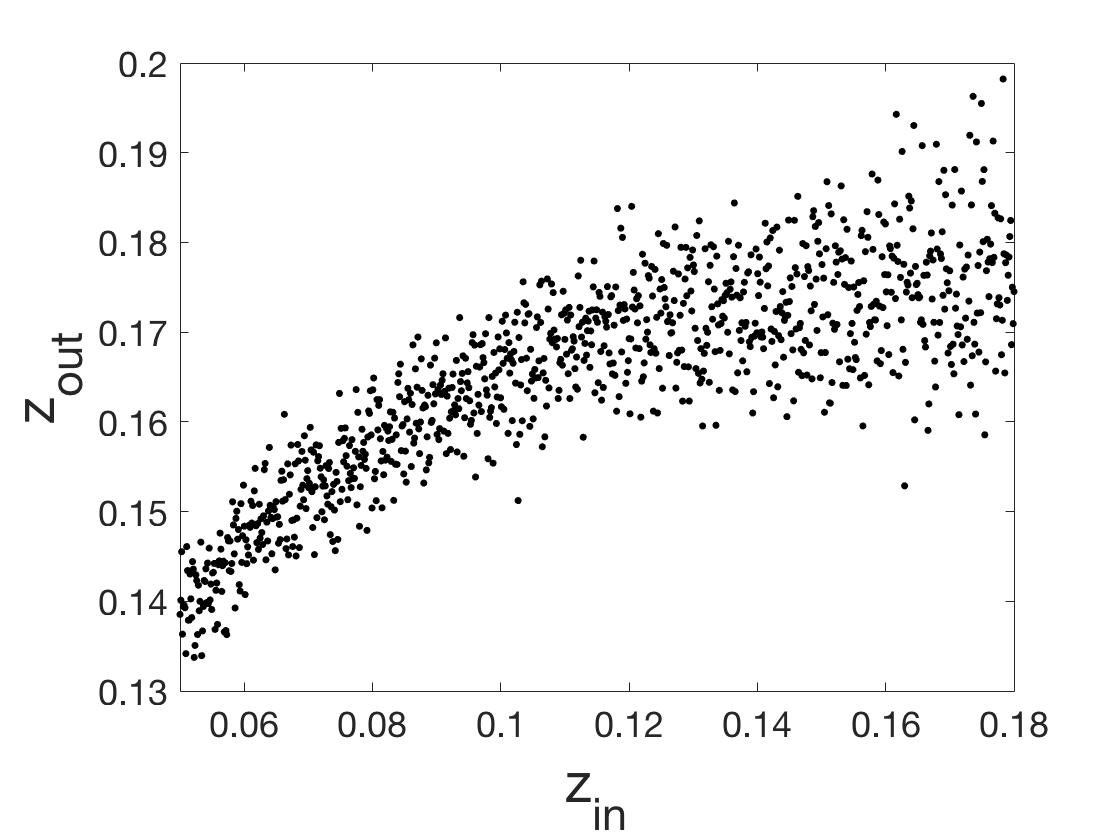} }}%
\caption{Left column: Image of intersections of a grid of $1000$ trajectories 
with $\Sigma$. 
 Right column: Final values of $z$ plotted against their 
initial values. 
}
    \label{return_maps_stochastic_1}%
\end{figure}

Next, we numerically study the effect of noise on the return map which was
 discussed in Section~\ref{sec:deterministic}.
We fix $h$ at $2.3$  and the other parameter values as in (\ref{par_values}). We then consider  
a line segment $\mathcal{L}$ on a plane $\Sigma$ far from $\mathcal{F}^{+}$,
integrate a sample path starting on equally spaced points on $\mathcal{L}$ forward, and record its first 
intersection with $\Sigma$. Fig.~\ref{return_maps_stochastic_1} shows the image of $\mathcal{L}$,
where we choose $\Sigma$ to be $x=0.18$ and 
$\mathcal{L}$ to be a line segment on $\Sigma$ with its $y$-coordinate fixed  at $0.22$,
and $z$ varying between $0.05$ and $0.18$. 
 We observe that different 
deterministically discernible regions (cf. Fig.~\ref{return_maps_stochastic_1}(b) 
and (f)), tend to look similar if more noise is added, i.e., the plateau disappears. 
Furthermore, the return maps look very regular, in certain parts even almost linear,
which relates to similar observations made in~\cite{K} for a simpler model.


\section{Analysis near the Folded Singularity}
\label{sec:fold_analysis}

The remaining step is to link our results to theoretical background, which describes
whether the SAOs are purely due to noise or have been induced by deterministic
mechanisms, e.g., by a folded node. Key theoretical steps to understand this
effect were made in~\cite{BGK}. In~\cite{BGK}, the authors study a slow-fast 
stochastic system in its normal form, where they quantify the effect of additive 
noise and obtain critical noise intensities beyond which the SAOs
 become hidden by fluctuations. Hence, we may ask whether we can 
  transform our system near a folded node into a suitable normal form SDE.

We will focus our analysis near a folded singularity of the desingularized 
system~\eqref{dsing} and follow ideas in~\cite{Wechselberger}, now extended to
the stochastic context. Let $p^{\star}=(x^{\star}, y^{\star}, z^{\star})\in 
\mathcal{F}^+$ be a fold point where the transversality condition fails, 
i.e. $(x^{\star}, y^{\star}, z^{\star})$ satisfy the following: 
\begin{eqnarray*}
  \left\{
\begin{array}{ll} 
  1-x^{\star}-\frac{y^{\star}}{\beta_1+x^{\star}}-
	\frac{z^{\star}}{\beta_2+x^{\star}} &=0\\
-1+\frac{y^{\star}}{(\beta_1+x^{\star})^2}+ 
\frac{z^{\star}}{(\beta_2+x^{\star})^2} &=0\\
\frac{y^{\star}}{\beta_1+x^{\star}}\left(\frac{x^{\star}}{\beta_1+
x^{\star}}-c\right)+\frac{z^{\star}}{\beta_2+x^{\star}} 
\left(\frac{x^{\star}}{\beta_2+x^{\star}}-d-hz^{\star}\right) &=0,
 \end{array} 
\right. 
\end{eqnarray*}
with $x^{\star}>c\beta_1/(1-c), y^{\star}, z^{\star}>0$. The fold curve 
$\mathcal{F}^+$ in its parametric form is given by
\bess \mathcal{F}^+ &=&
\left \{(s, \phi(s), \psi(s))  :  a\leq s\leq b \right\},  
\textnormal {where}\\
\phi(s) &=& \frac{(\beta_1+s)^2}{\beta_1-\beta_2}(1-\beta_2-2s), \\
\psi(s) &=& \frac{(\beta_2+s)^2}{\beta_1-\beta_2}\left(2s+\beta_1-1\right),\\
a &=& \frac{1-\max\{\beta_1, \beta_2\}}{2}, \ 
b =\frac{1-\min\{\beta_1,\beta_2\}}{2}.
\eess
Let $s^{\star}\in (a, b)$ be such that $s^{\star}=x^{\star}$, 
$\phi(s^{\star})=y^{\star}$ and $\psi(s^{\star})=z^{\star}$, 
i.e. $s^{\star}$ gets mapped to $p^{\star}$. Then for $z$ near 
$z^{\star}$, we can show that $\mathcal{F}^+$ can be linearly 
approximated by $\{(c_{11}z, c_{22}z, z)\}$, where
\begin{eqnarray*} c_{11}=\frac{\beta_1-\beta_2}{2(\beta_2+
x^{\star})(3x^{\star}+\beta_1+\beta_2-1)},\ c_{22}=-\frac{\beta_1
+x^{\star}}{\beta_2+x^{\star}},
\end{eqnarray*}
provided that $\beta_1\neq \beta_2$. 
Note that near $p^{\star}$, it is easy 
to see that $F_i$ can be approximated by the functions $G_i$ with 
$F_i^2(x,y,z)= G_i^2(x,y,z)+ O((x-x^{\star})^2, (y-y^{\star})^2, 
(z-z^{\star})^2)$  for $i=1,2,3$, where
\begin{eqnarray*}
G_1(X, Y, Z) &=& \sqrt{X+X^2+\frac{XY}{\beta_1+X}+\frac{XZ}{\beta_2+X}}, \\
G_2(X, Y, Z) &=& \sqrt{\frac{XY}{\beta_1+X}+c }, \\
G_3(X, Y, Z) &=& \sqrt{\frac{XZ}{\beta_2+X}+dZ+hZ^2 }.
\end{eqnarray*}
For the ease of calculations, we will consider system~\eqref{stoc} 
with $F_i$ replaced by $G_i$. Throughout our work for each positive 
integer $k$, $O(k)$ represents multivariate terms of the form 
$X^{\alpha}Y^{\beta} Z^{\gamma}$ such that $\alpha+\beta+\gamma=k$, 
where $\alpha, \beta, \gamma$ are non-negative integers.

\begin{proposition} 
Near the folded node singularity $(x^{\star}, y^{\star}, z^{\star})$ there 
exists a smooth change of coordinates such that system~\eqref{stoc} can be 
written in its normal form:
\begin{eqnarray*}   \left\{
\begin{array}{ll} \textnormal{d}X_t &= F_1(X_t, Y_t,Z_t) 
\textnormal{d}t+\kappa \sigma_1 
G_1(X_t, Y_t, Z_t)\ \textnormal{d}W_t^{(1)} \\
&+\sqrt{\zeta}\kappa \sigma_3  c_{11} G_3(X_t, Y_t, Z_t)\ 
\textnormal{d}W_t^{(3)} \\
\textnormal{d}Y_t &= \zeta F_2(X_t, Y_t,Z_t) \textnormal{d}t+
\frac{\kappa x^{\star}\sqrt{\zeta}}{\beta_1+x^{\star}}
\Big[(-\sigma_2  G_2(X_t, Y_t, Z_t) \\
&+O(1)) \textnormal{d}W_t^{(2)}+ (\sigma_3 c_{22}G_3(X_t, Y_t, Z_t)\\
&+O(1)) \textnormal{d}W_t^{(3)}\Big]\\
\textnormal{d}Z_t &= \zeta F_3(X_t, Y_t,Z_t) dt+\sqrt{\zeta}
\sigma_3 G_3(X_t, Y_t, Z_t)\ \textnormal{d}W_t^{(3)},
 \end{array} 
\right. \end{eqnarray*}
where the functions are given by
\begin{eqnarray*}  \begin{pmatrix}
F_1  \\
F_2\\
F_3
 \end{pmatrix}= \begin{pmatrix}
Y_t+X_t^2+CX_tY_t+\zeta O(X_t,Y_t,Z_t,\zeta)+\\
O(3) \\
A_1X_t+A_2Y_t+A_3Z_t+A_4 G_3^2(X_t, Y_t, Z_t)+\\
 O(2)\\
B_0+B_1X_t+B_2Z_t+O(2)
 \end{pmatrix},
\end{eqnarray*}

\begin{eqnarray*}   
\begin{pmatrix}
G_1^2  \\
G_2^2\\
G_3^2
 \end{pmatrix}= \begin{pmatrix}
2x^{\star}+\frac{2}{\kappa}X_t-Y_t+2c_{11}Z_t-\zeta C_0+O(2) \\
D_0+D_1X_t+D_2 Y_t+D_3Z_t-D_4\zeta +O(2) \\
E_0+\frac{\beta_2 z^{\star}}{\kappa(\beta_2+x^{\star})}X_t
+\left(\frac{\beta_2c_{11}z^{\star}}{(\beta_2+x^{\star})^2}+
hz^{\star} \right) Z_t +\\
O(2)
 \end{pmatrix},
\end{eqnarray*}
with $G_1$, $G_2$ and $G_3$ being the positive square roots of $G_1^2$, $G_2^2$ 
and $G_3^2$. The constants are all explicitly computable and given 
in~\eqref{a1}-\eqref{c22}.
\end{proposition}

\begin{figure}[h]%
    \centering
      \subfloat[$\chi_k$ vs $h$]
			{{\includegraphics[width=5.9cm]{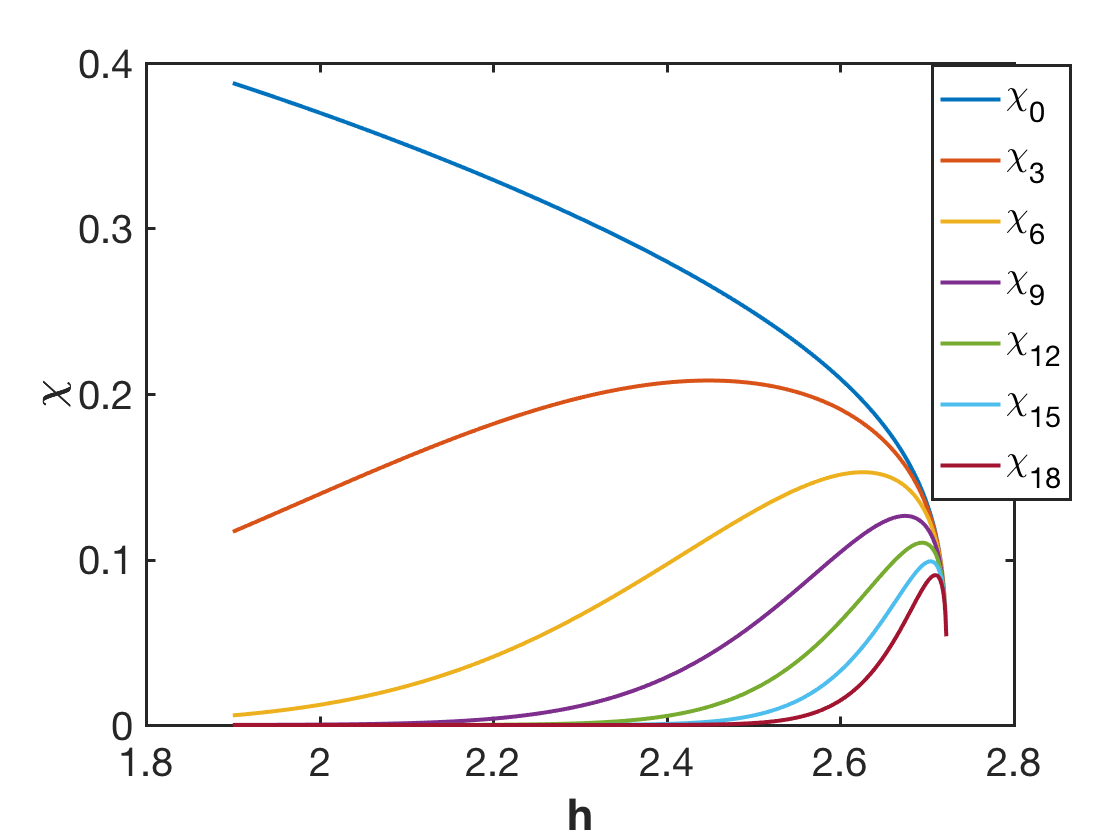} }}%
          \qquad
 \subfloat[A zoomed view] {{\includegraphics[width=5.9cm]
{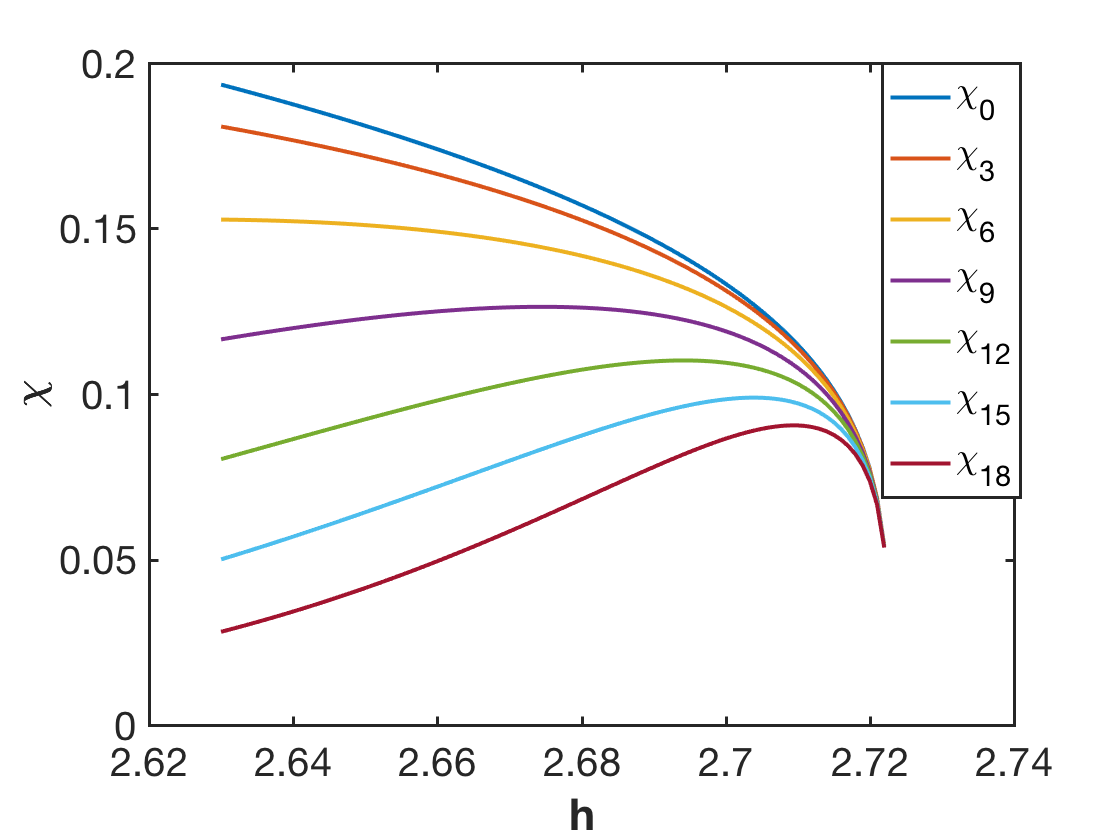}}}%
  \caption{Plots of $\chi_k$ as functions of $h$ for different values 
	of $k$ and parameter values given by (\ref{par_values}).}    
  \label{sig}%
\end{figure}

The proof is given in the Appendix. With this normal form, we have 
a very explicit tool to understand the parameter influence on the noise 
terms near the folded node. For instance, we consider the
parameter
\bess
c_{11}=\frac{\beta_1-\beta_2}{2(\beta_2+x^{\star})(3x^{\star}+\beta_1+\beta_2-1)}.
\eess 
If the two predators are similar with respect to their 
semi-saturation constants, i.e., $\beta_1\approx \beta_2$, then
fluctuations for the $Y$-predator near the folded node depend to leading-order 
only very weakly on the demographic noise of the $Z$-predator; similar ecological 
parameter studies can now be carried out for various parameters and SAO counts, 
but we leave this aspect for future work.  

In~\cite{BGK}, the noise level influence was carried out relative to 
the various canard orbits. Consider the $k$-th secondary 
canard. For each $k\in \mathbb{N}$, let 
\begin{equation}
\label{eq:BGKres}
\chi_k(\mu)=\mu^{1/4}\textnormal{e}^{(-(2k+1)^2 \mu)}
\end{equation}
where $\mu$ is the eigenvalue ratio at the folded node defined 
in~\eqref{eq:evratio} and $\chi_k$ is the noise level for additive
noise near a folded node in the normal form defined~\cite{BGK}. 
The curves $\chi_k(\mu)$ define the boundary regions 
in $(\mu,\chi)$-parameter space~\cite{BGK,BerglundGentzKuehn1}. To 
link this normal form result to our context, we plot $\mu=\mu(h)$
as a function of the main bifurcation parameter $h$ as shown in 
Fig.~\ref{sig}. Basically, upon decreasing the noise
level, within each region, we can distinguish more and more deterministic 
SAOs from pure-noise SAOs near a folded node. 

\begin{figure}[ht!]%
    \centering
			\includegraphics[width=6.0cm]{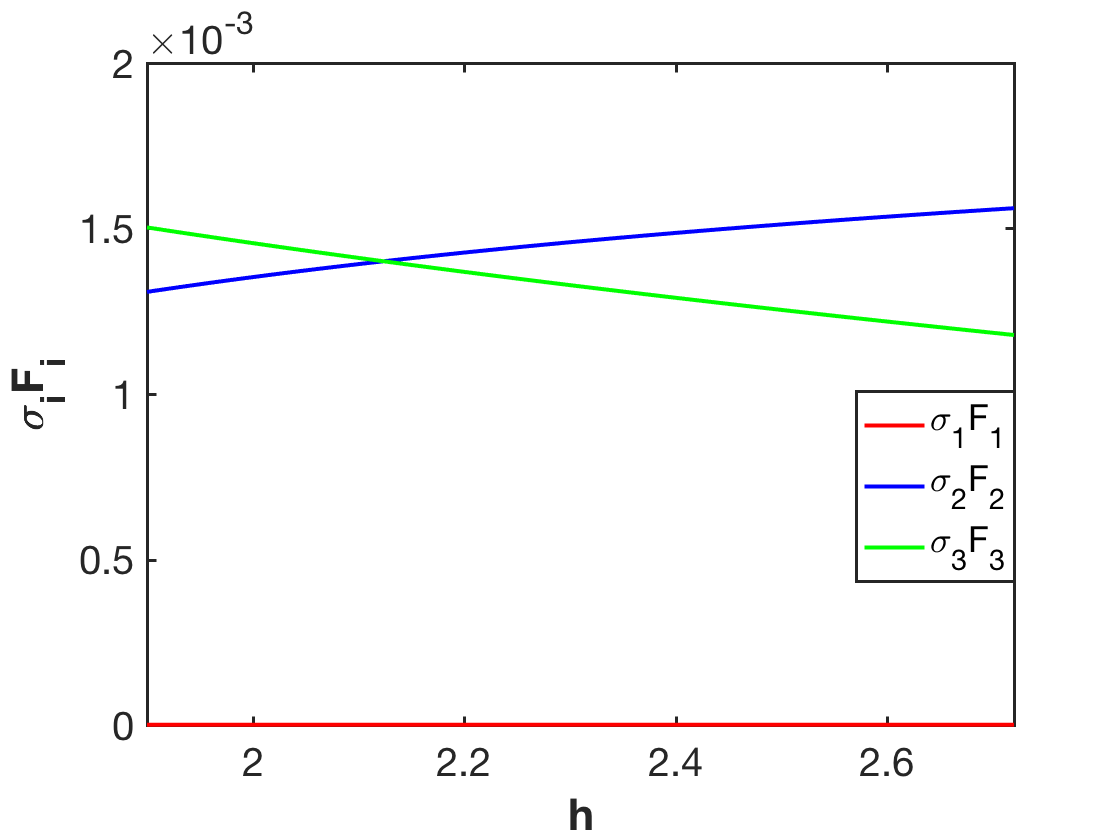}%
  \caption{The noise pre-factors $\sigma_iF_i$ corresponding to system (\ref{stocorig})
  evaluated at  $(x^{\star},y^{\star},z^{\star})$.
 }    
  \label{noise_prefactors}%
\end{figure}

To actually link the results completely, we have to determine 
the noise level $\chi$, which can be done by evaluating the products 
$\sigma_i F_i$ from~\eqref{stocorig} at the folded node. The results 
are shown in Fig.~\ref{noise_prefactors} for parameter values given by (\ref{par_values}).
 Our noise is rather small relative to 
$\chi$, so we may only expect occasional stochastic deviations from
SAOs in the regime far from the singular Hopf bifurcation.
However, the time series in Fig.~\ref{timeseries_0} show
more irregular behavior, which the theory does not cover
as the SAOs occur near a  singular Hopf bifurcation. The discrepancy is particularly apparent in 
the case when a small deterministically stable limit cycle exists; 
see Fig.~\ref{3d_samplepath}. 
Hence, we conclude that a singular Hopf bifurcation can still produce 
quite irregular MMOs, even for relatively small noise levels, where
a folded node normal form theory would not give an indication for the 
influence of noise on MMOs. To investigate this effect mathematically
remains a formidable challenge for future work.

\section{Appendix}
\label{sec:appnd}

We perform the following series of transformations to arrive at the normal 
form. We are going to need the following constants to simplify the notation:
\begin{eqnarray}
\label{a1} A_1 &=& -x^{\star}\left(\frac{\beta_1y^{\star}}
{(\beta_1+x^{\star})^3}+\frac{\beta_2z^{\star}}{(\beta_2+x^{\star})^3} 
\right), \\
 \label{a2} A_2 &=& \frac{\beta_1C_0 }{x^{\star}(\beta_1+x^{\star})}+C_1,\\
\label{c0} C_0 &=& c_{11}z^{\star} \left(\frac{x^{\star}}{\beta_1
+x^{\star}}-d-hz^{\star} \right),\\
  \label{c1} C_1 &=& \frac{x^{\star}}{\beta_1+x^{\star}}-c, \\
\nonumber A_3 &=&x^{\star}\kappa \Big[\frac{1}{\beta_2+x^{\star}}
\left(d+2hz^{\star}-c-\frac{(\beta_1-\beta_2)x^{\star}}{(\beta_1
+x^{\star})(\beta_2+x^{\star})}\right) \\
\label{a3} &-& c_{11}\left(\frac{\beta_1y^{\star}}{(\beta_1+x^{\star})^3}
+\frac{\beta_2z^{\star}}{(\beta_2+x^{\star})^3} \right) \Big],\\
\label{a4} A_4 &=& \frac{c_{11}c_{22}\kappa \beta_1\sigma_3^2 \zeta}
{(\beta_1+x^{\star})^2} ,\\
\label{b0} B_0 &=&z^{\star}\left(\frac{x^{\star}}{\beta_2+x^{\star}}
-d-hz^{\star}\right),\\
\label{b1} B_1 &=& \frac{\beta_2 z^{\star}}{(\beta_2+x^{\star})^2},\\
  \label{b2}  B_2 &=& \frac{x^{\star}}{\beta_2+x^{\star}}-d-2hz^{\star}
	+\frac{c_{11}\beta_2z^{\star}}{(\beta_2+x^{\star})^2} ,\\
\label{c} C &=& \frac{\beta_1}{\kappa x^{\star}(\beta_1+x^{\star})},\\
\label{d0}  D_0 &=& y^{\star}\left(\frac{x^{\star}}{\beta_1+x^{\star}}
+c \right),\\
\label{d1} D_1 &=&\frac{\beta_1y^{\star}}{\kappa(\beta_1+x^{\star})^2},\\
\label{d2} D_2 &=& -\frac{c\beta_1+(1+c)x^{\star}}{x^{\star}},\\
\label{d3}  D_3 &=& \frac{c_{11}\beta_1y^{\star}}{(\beta_1+x^{\star})^2}, \\
 \label{d4} D_4 &=& \frac{ c_{11} B_0}{x^{\star}} (c\beta_1+(1+c)x^{\star}),\\
  \label{e0}E_0 &=& z^{\star}\left(\frac{x^{\star}}{\beta_2+x^{\star}}+d
+hz^{\star} \right),\\
\label{k} \kappa &=& -x^{\star}\left(\frac{y^{\star}}{(\beta_1+x^{\star})^3}
+\frac{z^{\star}}{(\beta_2+x^{\star})^3} \right),\\
\label{c11} c_{11} &=& \frac{\beta_1-\beta_2}{2(\beta_2+x^{\star})(3x^{\star}
+\beta_1+\beta_2-1)},\\
 \label{c22} c_{22} &=&-\frac{\beta_1+x^{\star}}{\beta_2+x^{\star}}.
\end{eqnarray}
The main procedure consists of four main steps.

1) Translation to the origin: $X\mapsto X+x^{\star},\ Y\mapsto Y+y^{\star},
\ Z\mapsto Z+z^{\star}$. The transformed system reads as
\begin{eqnarray*}\label{stoc1}    
\left\{
\begin{array}{ll} \textnormal{d}X &=  f_1(X, Y, Z)\ \textnormal{d}t 
+\sigma_1 G_1(X,Y, Z)\ \textnormal{d}W_t^{(1)}\\
\textnormal{d}Y &= \zeta f_2(X,Y, Z)\ \textnormal{d}t
+\sigma_2\sqrt{\zeta}G_2(X, Y, Z)\ \textnormal{d}W_t^{(2)}\\
 \textnormal{d}Z &= \zeta f_3(X, Y, Z)\ \textnormal{d}t
+\sigma_3 \sqrt{\zeta}G_3(X, Y, Z)\ \textnormal{d}W_t^{(3)},
 \end{array} 
\right. 
\end{eqnarray*}
where the drift are specified further by the functions
\begin{eqnarray*}   \begin{pmatrix}
f_1  \\
f_2\\
f_3
 \end{pmatrix}= \begin{pmatrix}
(X+x^{\star})\Big(1-(X+x^{\star})-\frac{Y+y^{\star}}{\beta_1
+X+x^{\star}}\\
-\frac{Z+z^{\star}}{\beta_2+X+x^{\star}}\Big)\\
(Y+y^{\star})\left(\frac{X+x^{\star}}{\beta_1+X+x^{\star}}-c \right)\\
(Z+z^{\star})\left(\frac{X+x^{\star}}{\beta_2+X+x^{\star}}-
d-h(Z+z^{\star}) \right)
 \end{pmatrix},
\end{eqnarray*}
and the diffusion terms are fully defined via the functions 
\begin{eqnarray*} \begin{pmatrix}
G_1  \\
G_2\\
G_3
 \end{pmatrix}= \begin{pmatrix}
\sqrt{X+x^{\star}}\sqrt{1+(X+x^{\star})+\frac{Y+y^{\star}}{\beta_1
+X+x^{\star}}+\frac{Z+z^{\star}}{\beta_2+X+x^{\star}}} \\
\sqrt{Y+y^{\star}}\sqrt{\frac{X+x^{\star}}{\beta_1+X+x^{\star}}+c } \\
\sqrt{Z+z^{\star}}\sqrt{\frac{X+x^{\star}}{\beta_2+X+x^{\star}}
+d+h(Z+z^{\star}) }
 \end{pmatrix}.
\end{eqnarray*}

2) Rectification of the fold curve to the $z$-axis: 
$X=\tilde{X}+c_{11}\tilde{Z},\ Y=\tilde{Y}+c_{22}\tilde{Z}, 
Z=\tilde{Z}$, where $c_{11}$ and $c_{22}$ are defined 
by~\eqref{c11}-\eqref{c22}. The transformed system reads as
\begin{eqnarray*}\label{stoc2}    
	\left\{
\begin{array}{ll} \textnormal{d}\tilde{X} &= 
\tilde{ f_1}(\tilde{X}, \tilde{Y}, \tilde{Z})\ \textnormal{d}t 
+\sigma_1 \tilde{G_1}(\tilde{X},\tilde{Y}, \tilde{Z})\ \textnormal{d}W_t^{(1)} \\
&-\sqrt{\zeta}c_{11}\sigma_3\tilde{G_3}(\tilde{X},\tilde{Y},\tilde{Z}) 
\ \textnormal{d}W_t^{(3)} \\
\textnormal{d} \tilde{Y} &= \zeta \tilde{f_2}(\tilde{X},\tilde{Y}, 
\tilde{Z})\ \textnormal{d}t+ \sqrt{\zeta} (\sigma_2 
\tilde{G_2}(\tilde{X},\tilde{Y}, \tilde{Z})\ \textnormal{d}W_t^{(2)}\\
& -c_{22}\sigma_3 \tilde{G_3}(\tilde{X},\tilde{Y},\tilde{Z})
\textnormal{d}W_t^{(3)}) \\
 \textnormal{d}\tilde{Z} &= \zeta \tilde{f_3}(\tilde{X},\tilde{Y}, 
\tilde{Z})\ \textnormal{d}t+\sigma_3 \sqrt{\zeta} 
\tilde{G_3}(\tilde{X},\tilde{Y}, \tilde{Z})\ \textnormal{d}W_t^{(3)},
 \end{array} 
\right. \end{eqnarray*}
where 
\begin{eqnarray*}   \begin{pmatrix}
\tilde{f_1}  \\
\tilde{f_2}\\
\tilde{f_3}
 \end{pmatrix}= \begin{pmatrix}
-\frac{x^{\star}}{\beta_1+x^{\star}}\tilde{Y}-\frac{D_3}{y^{\star}}
\tilde{Y}\tilde{Z}+\kappa \tilde{X}^2-\frac{D_3}{c_{11}y^{\star}}
\tilde{X}\tilde{Y}\\
-\zeta c_{11}L
+O(3)\\
-\frac{A_1}{x^{\star}}(\beta_1+x^{\star})\tilde{X}+
\left(\frac{x^{\star}}{(\beta_1+x^{\star})^3}-c\right) 
\tilde{Y} + M \tilde{Z} \\
+O(2)\\
B_0 +B_1\tilde{X}+B_2 \tilde{Z}+O(2)
 \end{pmatrix},
\end{eqnarray*}
and we also have
\begin{eqnarray*}
L(\tilde{X},\tilde{Y}, \tilde{Z}) &=& B_0 +B_1\tilde{X}+B_2\tilde{Z} 
+O(2), \\
M &=&   -\frac{A_3(\beta_1+x^{\star})}{x^{\star}\kappa},\\
\end{eqnarray*}
with $A_1$ defined by (\ref{a1}), $A_3$ by (\ref{a3}), $B_i$ 
for  $i=0, 1,2$ by (\ref{b0})-(\ref{b2}), $D_3$ by (\ref{d3}) 
and $\kappa$ by (\ref{k}). The functions 
 $\tilde{G_i}$ for $i=1,2,3$ are positive square roots of 
$\tilde{G_i}^2$, where
\begin{eqnarray*} \begin{pmatrix}
\tilde{G_1}^2  \\
\tilde{G_2}^2\\
\tilde{G_3}^2
 \end{pmatrix}= \begin{pmatrix}
2x^{\star}+2\tilde{X}+\frac{x^{\star}}{\beta_1
+x^{\star}}\tilde{Y}+2c_{11}\tilde{Z} +O(2) \\
D_0+\frac{\beta_1y^{\star}}{(\beta_1+x^{\star})^2}
\tilde{X}+\left(\frac{x^{\star}}{\beta_1+x^{\star}}
+c\right)\tilde{Y}+D_3\tilde{Z}\\
+O(2) \\
E_0+B_1\tilde{X}+z^{\star} \left(\frac{\beta_2 c_{11}}{(\beta_2
+x^{\star})^2}+h \right)\tilde{Z}+O(2) 
 \end{pmatrix}
\end{eqnarray*}
with $D_0$ and $E_0$ defined by (\ref{d0}) and (\ref{e0}) respectively.

3) Stretch and translation of the coordinates $\tilde{X}$ and $\tilde{Y}$: 
\begin{eqnarray*} \hat{X}=\kappa \tilde{X}, \hat{Y}=-\frac{x^{\star} 
\kappa}{\beta_1+x^{\star}}\tilde{Y}-\zeta c_{11}\kappa E_0,\ \hat{Z}=\tilde{Z}.
\end{eqnarray*}
 The transformed system reads as
 \begin{eqnarray*}\label{stoc3}    
	\left\{
\begin{array}{ll} \textnormal{d}\hat{X} &= 
 \hat{f_1}(\hat{X}, \hat{Y}, \hat{Z})\ \textnormal{d}t +
\kappa \Big[\sigma_1 \hat{G_1}(\hat{X}, \hat{Y}, \hat{Z})\ 
\textnormal{d}W_t^{(1)}\\
&-\sqrt{\zeta}c_{11}\sigma_3\hat{G_3}(\hat{X}, \hat{Y}, 
\hat{Z}) \ \textnormal{d}W_t^{(3)}\Big] \\
\textnormal{d} \hat{Y} &= \zeta \hat{f_2}(\hat{X}, \hat{Y}, 
\hat{Z})\ \textnormal{d}t-\frac{x^{\star}\kappa 
\sqrt{\zeta}}{\beta_1+x^{\star}}\Big[\sigma_2 \hat{G_2}(\hat{X},
\hat{Y}, \hat{Z})\ \textnormal{d}W_t^{(2)}\\
& -c_{22}\sigma_3 \hat{G_3}(\hat{X}, \hat{Y}, 
\hat{Z})\textnormal{d}W_t^{(3)}\Big] \\
 \textnormal{d}\hat{Z} &= \zeta \hat{f_3}(\hat{X}, \hat{Y},
 \hat{Z})\ \textnormal{d}t+\sigma_3 \sqrt{\zeta}\hat{G_3}(\hat{X}, 
\hat{Y}, \hat{Z})\ \textnormal{d}W_t^{(3)},
 \end{array} 
\right. \end{eqnarray*}
where 
\begin{eqnarray*}   \begin{pmatrix}
\hat{f_1}  \\
\hat{f_2}\\
\hat{f_3}
 \end{pmatrix}= \begin{pmatrix}
\left(1+c_{11}\kappa C \hat{Z} \right)\hat{Y} 
+\hat{X}^2 +\kappa C \hat{X}\hat{Y}+ O(\zeta)+O(3)\\
A_1 \hat{X}+ C_1\hat{Y}+A_3 \hat{Z} +\zeta c_{11}\kappa B_0 C_1+O(2)\\
B_0 +B_1\hat{X}+B_2\hat{Z}+O(2)
 \end{pmatrix},
\end{eqnarray*}
and 
\begin{eqnarray*} \begin{pmatrix}
\hat{G_1}^2  \\
\hat{G_2}^2\\
\hat{G_3}^2
 \end{pmatrix}= \begin{pmatrix}
2x^{\star}+\frac{2}{\kappa}\hat{X}-\frac{1}{\kappa}\hat{Y}
-\zeta c_{11}B_0 +2c_{11}\hat{Z} +O(2) \\
D_0+D_1\hat{X}+D_2\hat{Y}+D_3\hat{Z}+\zeta c_{11}B_0 D_2\\
+O(2) \\
B_0+\frac{B_1}{\kappa}\hat{X}+\left(\frac{\beta_2 
c_{11}z^{\star}}{(\beta_2+x^{\star})^2}+h z^{\star} \right)\hat{Z}+O(2) 
 \end{pmatrix}
\end{eqnarray*}
with $C$ defined by (\ref{c}), $D_i$ by (\ref{d0})-(\ref{d3}) 
and $\hat{G_i}$ being the positive square root of $\hat{G_i}^2$ 
for $i=0,1,2,3$.

4) A final transformation of $\hat{X}$ and $\hat{Y}$:  
$\bar{X}=\hat{X},\ \bar{Y}=\left(1+\frac{c_{11}\beta_1}
{x^{\star}(\beta_1+x^{\star})}\hat{Z}\right) \hat{Y},\ \bar{Z}=\hat{Z}$. 
Carrying out this transformation, dropping the overbars, 
and denoting the components by 
$(X_t, Y_t, Z_t)$, we thus obtain the proposition.

\end{document}